\documentclass{svjour3}
\usepackage{color}
\smartqed  
\usepackage{epstopdf}
\usepackage{subfigure}

\usepackage[english]{babel}
\usepackage{graphicx}
\usepackage{mathrsfs}
\usepackage{amsfonts}
\usepackage{fancyhdr}
\usepackage{textcomp}
\usepackage[official]{eurosym}
\usepackage{epstopdf}
\usepackage{color,soul}
\usepackage{mathdots}

\newcommand\NN{{\hbox{I\kern-.14em{N}}}}
\newcommand\RR{{\hbox{I\kern-.14em{R}}}}
\newcommand\ZZ{{\hbox{I\kern-.14em{Z}}}}

\begin{document}

\title{A scenario-based framework for supply planning under uncertainty: stochastic programming versus robust optimization  approaches 
}

\titlerunning{Supply planning under uncertainty}        

\author{Francesca Maggioni        \and
        Florian A. Potra \and
				Marida Bertocchi   
}


\institute{F. Maggioni \at
              Department of Management, Economics and Quantitative Methods,
University of Bergamo, Bergamo, Italy \\
              Tel.: +39-0352052649\\
              Fax: +39-0352052549\\
              \email{francesca.maggioni@unibg.it}           
           \and
           F. A. Potra \at
              Department of Mathematics \& Statistics, University of Maryland, Baltimore County, U.S.A.
							\and
							M. Bertocchi \at  Department of Management, Economics and Quantitative Methods,
University of Bergamo, Bergamo, Italy
}

\date{Received: date / Accepted: date}

\maketitle

\begin{abstract}
In this paper we analyze the effect of two modelling approaches for supply planning problems under uncertainty: two-stage stochastic programming (SP) and robust optimization (RO).
The comparison between the two approaches is performed through a  \textit{scenario-based framework} methodology, which  can be applied to any optimization problem affected by uncertainty.
For SP we compute the minimum expected cost based on the specific probability distribution of the uncertain parameters related to a  set of  scenarios.
For RO we consider  static approaches where random parameters belong to box or ellipsoidal uncertainty sets in compliance with the data used to generate SP scenarios.
Dynamic approaches for RO, via the concept of  adjustable robust counterpart, are also considered.

The efficiency of the methodology  has been illustrated for a supply planning problem  to optimize vehicle-renting and procurement transportation activities involving uncertainty on demands and on buying costs for extra-vehicles.
Numerical experiments through the scenario-based framework allow a fair comparison in real case instances. Advantages and disadvantages of RO and SP are discussed.

\end{abstract}

\section{Introduction} 

In this paper we analyze the effect of two modelling approaches, Stochastic Programming (SP) and Robust Optimization (RO), to a supply planning problem under uncertainty.
Stochastic Programming  and Robust Optimization  are considered two alternative techniques to deal with uncertain data both in a single period and  in a multi-period decision making process.
The main difficulty associated with the  former is the need to provide the probability distribution functions of the underlying stochastic parameters.
This requirement creates a heavy burden to the user because in many real-world  situations, such information is unavailable or hard to obtain (see for example \cite{BL11} and  \cite{RS03}).
On the other side RO  addresses the uncertain nature of the problem without making specific assumptions on probability distributions: the uncertain parameters are assumed to belong to a deterministic uncertainty set.  The drawback of this approach is the potentially strong dependence of the solution on the rather arbitrarily chosen uncertainty set.
RO adopts a min-max approach that addresses uncertainty by guaranteeing the feasibility and optimality of the solution against all instances of the parameters within the uncertainty set \cite{BBC2011,D98}.
A vast literature about the hypotheses that have to be imposed on the structure of the uncertainty set in order to have computationally tractable problems are available, see \cite{BS04}  and \cite{So73}   for polyhedral uncertainty sets and  \cite{BTN99,BGGN04,EL97,EOL98} for ellipsoidal uncertainty sets.
The original RO model deals with static problems where all the decision variables have to be determined before any of the uncertain parameters are realized.
This is not the typical situation in most problems that are multi-period in nature, and where a decision at any period can and should account for data realizations in previous periods.
An extension of robust optimization to a dynamic framework was analyzed in  \cite{BTN99,BGGN04,BTEGN09,BG12,BB2010,BG2010} and many others via the concept of adjustable robust counterpart (ARC) and affinely adjustable robust counterpart (AARC),  where  part of the decision variables, the so-called adjustable variables, have to be determined after a portion of the uncertain data is realized. In the case of AARC the dependence of  the adjustable variables on the realized data is represented by an  affine function. The introduction of AARC is motivated by the fact that in most of the cases the ARC approach is computationally intractable. 

Comparing SP and RO methods considering a cost-based approach is
somehow unilateral. The philosophy of the two
approaches to deal with uncertainty is  completely different indeed. The
RO approach is well-suited to the cases where the
optimizer wants to hedge the result against all imaginable outcomes
of the uncertain events. Such
methods can be successfully implemented usually only in the cases
where the hedging against rare events is cheap or at
least at low cost with respect to possible consequences. On the
other side, the SP approach includes the probabilistic information about the events into consideration. As a
consequence, less probable events are counted only with small
weights implying less conservative and, in many cases, cheaper
solutions. As drawback, we have a risk of appearance of rare
but expensive events which can result in high actual costs if such
an ``unhappy'' event is realized.

To  make a fair comparison between the above approaches we consider a \textit{scenario-based framework},
which does not just compare the different solution approaches in terms of costs;
it compares the costs of implementing the non-adjustable (or first-stage) solutions obtained using the information available up to now, with SP or RO.
The adjustable variables are then determined when the uncertainty is revealed by solving a deterministic problem.
This heuristic methodology allows also to overcome the typical problem of computational intractability of ARC and can be applied to any optimization problem affected by uncertainty.



The efficiency of the methodology has been illustrated  on a supply planning problem  to optimize vehicle-renting and procurement transportation activities to satisfy demand in several destinations out of several origins.
Uncertainty on demands and  cost of extra vehicles is considered.
The problem consists in  determining  the number of vehicles to book, at the end of each time period, from each plant  of the  set of suppliers, to replenish 	a certain good at factories   in order to minimize the total cost, given by the sum of the transportation costs from origin to destinations (including the discount  for vehicles booked but not used) and the cost of buying  units of product from external sources in case of inventory shortage at the destinations. We formulate and solve the problem as  a two-stage stochastic programming model and as  robust optimization models with different uncertainty sets where scenarios of demand and buying costs are built on historical data.

Since the demand of goods is in general highly affected by the economic conditions, a reliable forecast and reasonable estimates of probability distributions are difficult to obtain.
Furthermore, the supply-planning company would avoid to negotiate the quantity of vehicles with the suppliers every time period,  being immunized against every possible demand realization allowing to save a lot of operational activities.
This motivates us  to consider besides the SP approach, also  RO approaches  and to quantify the value (or extra-cost) of RO guarantees.
First we consider  static approaches where the uncertain parameters belong to box or box-ellipsoidal uncertainty sets, and then dynamic approaches, via the concept of ARC with scenario generated uncertainty set and scenario-based framework methodology. A robust solution at the tactical level allows to find a feasible solution for the operational planning problem for each possible realization of demand in the uncertainty set considered.

Both SP and RO allow  to determine the nonadjustable variables, i.e., the number of vehicles to book at the end of each time period, using the information available at that time. As new information on demand and buying costs from external sources become available, the adjustable (or recourse) decision variables have to be determined.
We describe five strategies for updating the adjustable variables given the values of the already determined values of the nonadjustable variables. The methodology allows to quantify  the cost saving of the SP approach compared to the RO, as well as  the value of a more conservative strategy avoiding a negotiation  of the number of vehicles every period with the suppliers or third-party service providers (3PL).  Numerical experiments through the scenario-based framework allow a fair comparison in real instances inspired by a gypsum replenishment problem of an Italian cement factory producer.

The paper is organized as follows:  Section \ref{sec:StochasticVersusRobustOptimization} introduces basic concepts about stochastic and robust optimization.
A scenario-based framework used to compare the two approaches is presented in Section \ref{sec:AScenarioBasedFrameworkForComparisonOfModellingApproachesUnderUncertainty} while
 Section \ref{SPUU} describes the supply planning problem.
Section \ref{Models} discusses the two-stage stochastic programming formulation, robust formulations and the scenario based framework for comparison.
Finally, Section \ref{sec:numericalresults} discusses the numerical results. Conclusions follow.

\section{Stochastic versus robust optimization: basic facts and notations}
\label{sec:StochasticVersusRobustOptimization}

We consider the uncertain linear optimization problem

\begin{equation}
\left\{\min_{x} \left\{c^{T} x + o:\  Ax \leq h\right\}\right\}_{(c,o,A,h)\in\mathcal{U}}\ ,
\label{RO1}
\end{equation}
as a collection of linear optimization problems with data varying in a given uncertainty set $\mathcal{U}$ where  $x \in\mathbb{R}^{n,+}$  is the vector of non-negative decision variables, $c \in\mathbb{R}^n$  and $o\in\mathbb{R}$ are the coefficients of the objective function, $A\in\mathbb{R}^{m \times n}$ is the constraint matrix, and $h\in \mathbb{R}^{m}$ is the right hand side vector.

Let us introduce a scenario based  stochastic programming formulation of problem (\ref{RO1}) which includes equality constraints, by rewriting the constraint  $Ax \leq h$ in (\ref{RO1}) as
$Ax + \delta  = h,
\label{equalityconstraint}
$
where  $\delta$ is a slack variable, or equivalently
$
\left[\begin{array}{c|c}
A  & I
\end{array}\right]
\left[\begin{array}{c}
x \\
\delta
\end{array}\right]=h,
$
where $I\in\mathbb{R}^{m\times m}$ is  the identity matrix.

Let us consider:
\begin{itemize}
	\item[-]  A finite set $\zeta^s$, $s\in\mathscr{S}=\left\{1,\dots,S\right\}$ of realizations (or scenarios) of a random event  $\zeta$ which affects the data varying in $\mathcal{U}$.
	The random process $\zeta$ is defined on a probability space $(\Xi,\mathscr{F},P)$  with support $\Xi$ and given probability distribution $P$ on the $\sigma-$algebra $\mathscr{F}$.

\item[-] A partitioning of the  decision variable $x=[u ; v^1;\ldots;  v^S]$,
where $u\in\mathbb{R}^{n_1,+}$ represents the first-stage decision  which has to be taken without full information on the
random event $\zeta$. When full information is received on the realization
of the random vector, then, second-stage or recourse actions $v^s\in\mathbb{R}^{n_2,+}$, $s=1,\dots,S$ are taken.
Throughout this paper, we use ``;'' for adjoining elements  in a column.
\item[-] A partitioning of vector cost $c=[p;P^1 q^1;\ldots;P^S q^S]$, 
where $p\in\mathbb{R}^{n_1}$ is the first stage cost, $q^s\in\mathbb{R}^{n_2}$, $s=1,\dots,S$, is the second stage cost and  $P^s$ is the probability of scenario $s=1,\dots,S$.
\item[-] A partitioning of matrix $\hat{A}=[A | I]$ as follows:
\begin{displaymath}
\hat{A}=\left[\begin{array}{c|ccc}
\tilde{A} & & \mathbf{0} &\\
\hline
T^1 & W & &\\
\vdots &  & \ddots & \\
T^S & & & W
\end{array}\right]\ ,
\end{displaymath}
where  $\tilde{A}\in\mathbb{R}^{m_1 \times n_1}$ is a deterministic matrix, $T^s \in \mathbb{R}^{m_2\times n_1}$, $s=1,\dots,S$ are called technology matrices, $W\in\mathbb{R}^{m_2 \times n_2}$ is the matrix of fixed recourse,  and $\mathbf{0}$ the null matrix of dimension $m_1 \times S n_2$.
\item[-]  The partitioning of right hand side vector $h=[\tilde{h};h^1;\ldots;h^S]$,
where $\tilde{h}\in\mathbb{R}^{m_1}$ is a deterministic right hand side vector and $h^s\in \mathbb{R}^{m_2}$, $s=1,\dots,S$ are  stochastic right hand side terms.
\end{itemize}
\begin{definition}
 The two-stage linear stochastic problem with fixed recourse, is formulated as follows:
\begin{eqnarray}
\label{sp_multi_scenario}
& &\! \min_{u,v^1,\dots,v^S} p^{T} u + \sum_{s=1}^{S}P^s \left({q^s}^{T} v^s  \right) + o \nonumber \\
& &\textrm{s.t. }	 \tilde{A} u=\tilde{h}\ ,\nonumber\\
& & \qquad T^s u + W v^s= h^s\ ,\quad s=1,\dots, S\ ,\\
& &\qquad u \geq 0\ ,\quad v^{s}\geq 0\ ,\quad s=1,\dots,S\ .\nonumber
\end{eqnarray}
\end{definition}
We note that $o$ usually does not appear in SP problems, however it is reported  for the sake of consistency of notation with the RO formulation.

Let us introduce the robust optimization counterpart of problem (\ref{RO1}).
We assume that the uncertainty set $\mathcal{U}$ is parametrized in an affine way by the perturbation vector $\zeta=[\zeta_1; \dots;
\zeta_L]$  in a given perturbation set $Z$:

\begin{equation}
\mathcal{U}=\left\{[c;o;A;h]=[\bar{c};\bar{o};\bar{A};\bar{h}] + \sum_{\ell=1}^L\zeta_\ell[c_\ell;o_\ell;A_\ell;h_\ell]: \zeta \in Z \subset\mathbb{R}^{L} \right\}\ ,
\end{equation}
where $[c_\ell;o_\ell;A_\ell;h_\ell]$ represent possible perturbations from the nominal data  $[\bar{c};\bar{o};\bar{A};\bar{h}]$.
\begin{definition}
The Robust Counterpart of the uncertain linear optimization problem (\ref{RO1}) is formulated as follows
\begin{equation}
\min_{x}\left\{ \sup_{(c,o,A,h)\in\mathcal{U}}c^{T} x + o:\  Ax \leq h,\ (c,o,A,h)\in\mathcal{U}\right\}\ ,
\label{RO2}
\end{equation}
that is minimizing the worst total cost over all feasible solutions.
\end{definition}

Problem (\ref{RO2}) can be equivalently formulated as
\begin{equation}
\min_{x,w}\left\{w:\  c^{T} x -w\leq - o,\  Ax \leq h,\  \forall (c,o,A,h)\in\mathcal{U}\right\}\ .
\label{RO3}
\end{equation}

Notice that if $x$ is a feasible solution of (\ref{RO2}), then $x$ remains feasible when we extend the uncertainty set $\mathcal{U}$ to its convex hull $Conv(\mathcal{U})$.


%

We consider tractable formulations of the robust problem (\ref{RO3}) with uncertainty set $\mathcal{U}$ computationally tractable, that is solvable in  polynomial time with a specific description of the uncertainty set. There are three well known formulations of RO problems in literature; these are given by    \cite{BTN00,BTEGN09,BS04} and \cite{So73}.   They all share the advantage that minimal assumptions about the nature of the uncertainties have to be made and they differ in the ways the uncertainty sets are represented. More specifically, the formulations by Soyster \cite{So73} and by Bertsimas and Sim \cite{BS04} use polyhedral uncertainty sets, while the formulation by Ben-Tal and Nemirovski
\cite{BTN99,BGGN04,EL97,EOL98}  considers an ellipsoidal uncertainty set, transforming the original LP problem into a Second Order Cone Programming (SOCP) problem.

Among possible  formulations we consider the case of box uncertainty sets and ellipsoidal uncertainty sets. To show the construction of the uncertainty sets we focus on a single uncertain linear inequality:
\begin{equation}
\left\{ a^{T} x \leq h \right\}_{(a,h)\in\mathcal{U}}\ .
\label{RO6}
\end{equation}
The tractable formulation of $\mathcal{U}$ through a box  uncertainty set is given by
\begin{equation}
a^{T} x\le h,\,\forall [a;h]\!\in\!\left\{\! [\bar{a};\bar{h}]\!+\!\sum_{\ell=1}^L\zeta_\ell [a^\ell;h^\ell]:\, \zeta\!=\![\zeta_1;\dots;\zeta_L]\!\in\!\mathbb{R}^L,\, \|\zeta\|_{\infty}\le 1\!\right\},
\label{RObox}
\end{equation}
where $\|\cdot\|_{\infty}$ is the infinity norm, while the formulation through an ellipsoidal uncertainty set is given by
\begin{equation}
a^{T} x\le h,\;\forall [a;h]\in\left\{ [\bar{a};\bar{h}]+\sum_{\ell=1}^L\zeta_\ell [a^\ell;h^\ell]\;:\; \zeta=[\zeta_1;...;\zeta_L]\in\mathbb{R}^L,\, \|\zeta\|_{2}\le \Omega\right\}\ ,
\label{ROellipsoid}
\end{equation}
where $\|\cdot\|_{2}$ is the Euclidean norm and $\Omega$ the radius.

As shown in \cite{BTEGN09} (\ref{ROellipsoid}) can be written as
\begin{equation}
\label{socc}
\left[(\bar{h}-\bar{a}^{T} x )/\Omega; h^1  - {a^{1}}^{T} x;\dots;h^L - {a^{L}}^{T} x \right]\in
\mathscr{LC}_{L+1},
\end{equation}
where $\mathscr{LC}_{L+1} =\left\{t= \left[t_0;t_1;\dots;t_L\right]\in\mathbb{R}^{L+1}: t_0 \geq \left\|t_1;\dots;t_L\right\|_2\right\}\ ,$
is  the second order (or Lorentz) cone of $\mathbb{R}^{L+1}$, so that that (\ref{ROellipsoid}) is equivalent to
\begin{equation}
\Omega \sqrt{\sum_{\ell=1}^L([a^\ell]^{T} x - h^\ell)^2}\leq \bar{h} - \bar{a}^{T} x\ .
\end{equation}

Let us consider a random vector $[a;h]$ defined by
\begin{equation}
[a;h]=[\bar{a};\bar{h}]+\sum_{\ell=1}^L\zeta_\ell [a^\ell;h^\ell],\label{rand1}
\end{equation}
where
\begin{equation}
\zeta_1,\!\ldots,\!\zeta_L: \mbox{zero mean independent random variables with values in } [-1,1]\ .\label{rand2}
\end{equation}
The following result holds true (see Corollary 2.3.2. in \cite{BTEGN09}):
\begin{proposition}\label{prop}
If $[a;h]$ is the random vector given by (\ref{rand1})-(\ref{rand2}) and $x$ is a solution of (\ref{socc}) then
\begin{equation}\label{prob}
\mbox{\rm Prob}\left\{\ a^Tx > h\right\}\le e^{-\frac{\Omega^2}{2}}.
\end{equation}
\end{proposition}
We note that the above result holds for any probability distribution for the random vector $\zeta=[\zeta_1;\ldots;\zeta_L]$ that satisfies (\ref{rand2}). Let us denote by ${\cal P}$ the family of all probability distributions that satisfy (\ref{rand2}) and consider the ambiguous chance constraint  where $\varepsilon\in (0; 1)$ is a prespecified small tolerance:
\begin{equation}\label{chance}
\forall P\in {\cal P}\quad {\mbox{\rm Prob}}_{\zeta\sim P}\left\{\zeta\; : \;
\bar{a}^{T} x+\sum_{\ell=1}^L\zeta_\ell [a^\ell]^{T} x\,>\,\bar{h}+\sum_{\ell=1}^L\zeta_\ell h^\ell\right\}\le\varepsilon .
\end{equation}
We note that the relation between $\epsilon$ and $\Omega$ is
$\Omega= \sqrt{\ln \epsilon^{-2}}. $

The above problem is called an {\it ambiguous} chance constraint  because we do not have any knowledge about the probability distribution $P$ except the fact that it belongs to the class ${\cal P}$ and we have uncertainty in the particular realization of the data (given its distribution). From Proposition~\ref{prop} it follows that any $x$ satisfying the second order cone constraint (\ref{socc}) is a solution of ({\ref{chance}).

All decision variables described above in case of RO approaches represent \textit{here and now decisions}, i.e., decisions that are taken before the actual data reveals itself. This is too restrictive since in reality there may be situations in which the decisions must adjust themselves to the actual data.  To model adjustability of the variables we can proceed as follows:
$\forall \, j\leq n$ we suppose that  $x_j$ is dependent on a portion of the data $(c,o,A,h)\in\mathcal{U}$, i.e.
\begin{equation}
x_j=X_j(\pi_j(\zeta)),\quad j=1,\dots,n,
\end{equation}
where  the uncertain data are denoted by $\zeta=(c,o,A,h)$, $\pi_j$	describes a portion of the uncertain data  by a suitable linear mapping on $\zeta$, and $X_j$ are decision rules to be chosen. We can now replace problem (\ref{RO3}) with the usage of the new decision rules $X_j$ as follows:
\begin{equation}
\min_{\left\{X_j(\cdot)\right\}_{j=1}^n,w}\left\{w:\  c(\zeta)^{T} X(\zeta) -w\leq - o(\zeta),\  A(\zeta) X(\zeta) \leq h(\zeta),\,  \zeta \in \mathcal{U}\right\}\ ,
\label{RO7}
\end{equation}
where $X(\zeta)=[X_1(\pi_1(\zeta));\dots;X_n(\pi_n(\zeta))]$.
The robust optimization problem (\ref{RO7})  is called Adjustable Robust Counterpart (ARC).

We consider:
\begin{itemize}
	\item the case of fixed recourse, i.e. for every adjustable variable  all its coefficients in the objective function and the left hand side of the constraints are certain;
	\item a scenario-generated uncertainty set described as a convex hull of finitely many scenarios $\zeta^s$, $s=1,\dots,S$.
\end{itemize}

Assuming  that $x=\left[u; v(\zeta)\right]$, 
where $u$ refers to the here-and-now variables (non-adjustable) and $v$ refers to the wait-and-see (adjustable) variables, (\ref{RO7}) becomes
\begin{equation}
\min_{u,v(\zeta),w}\!\!\left\{\!w\!:\! p(\zeta)^{T}u\! +\! q^{T}\! v(\zeta) \!-\! w \! \leq\! - o(\zeta),\   T(\zeta) u \!+\! W \! v(\zeta) \leq \! h(\zeta),\   \zeta \!\in \! Conv\left\{\zeta^1,\dots,\zeta^S\!\right\}\!\right\}
\label{RO8}
\end{equation}
where  $p(\zeta), T(\zeta),o(\zeta),h(\zeta)$ are affine in $\zeta$.

The following theorem holds (see \cite{BTEGN09}):
\begin{theorem}
Under the assumption of fixed recourse and  scenario-generated uncertainty set described above, the ARC (\ref{RO8}) is equivalent  to the computationally tractable problem
\begin{equation}
\min_{u,\left\{v^s\right\}_{s=1}^S,w}\hspace{1pt}\left\{\!w: \! p(\zeta^s)^{T}u\! + \!q^{T} v^s\! - \!w \!\leq \!- o(\zeta^s),\!  T(\zeta^s) u \!+\! W v^s \leq h(\zeta^s),\;   s=1,\!\dots\!,S\!\right\}\qquad
\label{RO9}
\end{equation}
and their optimal values are equal.
Moreover, if $\bar{w}$, $\bar{u}$, $\left\{\bar{v}^s\right\}_{s=1}^{S}$ is a feasible solution to (\ref{RO9}), then the pair $\bar{w}$, $\bar{u}$ augmented by the decision rule for the adjustable variables
\begin{equation}
v(\zeta)= \sum_{s=1}^{S}\lambda_s(\zeta)\bar{v}^s\ ,
\label{v}
\end{equation}
form a feasible solution to the ARC.
Here $\lambda(\zeta)=[\lambda_1(\zeta);\dots;\lambda_S(\zeta)]$ is a nonnegative vector with the unit sum of entries such that
\begin{equation}
\zeta=\sum_{s=1}^{S}\lambda_s(\zeta)\zeta^s.
\label{eqlambda}
\end{equation}
\label{Theorem1}
\end{theorem}
The assumption of fixed recourse is essential, since without it the ARC may become intractable (see \cite{BGGN04}).
However, from a practical point of view, equation (\ref{eqlambda}) is not easily satisfied unless the number of scenarios $S$ becomes extremely large.

\section{A scenario based framework for comparison of modelling approaches under uncertainty}
\label{sec:AScenarioBasedFrameworkForComparisonOfModellingApproachesUnderUncertainty}

In this section we propose a scenario based framework for comparison of stochastic programming and robust optimization approaches.  The framework can be applied to any optimization problem affected by uncertainty.


We introduce the following notation:
\begin{enumerate}
\item[{$\mathbf{M_1}$}] is the stochastic optimization problem (\ref{sp_multi_scenario});
\item[{$\mathbf{M_2}$}] is the robust optimization problem (\ref{RO3}) with box constraints as in (\ref{RObox});
\item[{$\mathbf{M_3}$}] is the robust optimization problem (\ref{RO3}) with ellipsoidal constraints as in (\ref{ROellipsoid});
\item[{$\mathbf{M_4}$}] is the computationally tractable robust optimization problem (\ref{RO9}).
\end{enumerate}

Let
\begin{equation}\label{barb}
\hat{\Delta} = \left\{\hat{\zeta}^1,\hat{\zeta}^2,\ldots,\hat{\zeta}^{S}\right\}\ ,
\end{equation}
be a given  set of  scenarios of the uncertain parameters $\zeta$ eventually obtained by historical data. We consider  a set of indices
\begin{equation}\label{bs}
\bar{\mathscr{S}}=\left\{1,\ldots,\bar{S}\right\}\subset\mathscr{S}\ ,
\end{equation}
with cardinality $\bar{S}<S$.
For each such $\tau=\bar{S},\dots,S-1$ we compute the quantities
\begin{equation}\label{quantd}
\bar{\zeta}^{\tau}=\frac{1}{\tau}\sum_{s=1}^{\tau} \hat{\zeta}^s,\qquad \textrm{ and } \qquad \max_{s=1,\dots,\tau}\left\|\hat{\zeta}^{s} - \bar{\zeta}^{\tau}\right\|_2 \ ,
\end{equation}
which represent respectively the nominal and variable terms in (\ref{RObox})  and (\ref{ROellipsoid}).

For each $\tau=\bar{S},\ldots,S-1$  we  find the optimal first stage, or nonadjustable decision variables  solutions $u^*$, of the corresponding optimization problem $\mathbf{M_m},\;\mathbf{m}=\mathbf{1},\ldots,\mathbf{4}$, using only the information contained in the vectors
\begin{equation}
\hat{\zeta}^1,\hat{\zeta}^2,\ldots,\hat{\zeta}^{\tau}.
\label{infobis}
\end{equation}
Assume now that the vectors $\hat{\zeta}^{\tau+1}$ become available.
Then we can solve the optimization problem
\begin{equation}
\min_{w,v}\left\{w:\  {p(\hat{\zeta}^{\tau+1})}^{T}u^{\ast}\! +\! q^{T} v \!- \!w \!\leq \!- o(\hat{\zeta}^{\tau+1}),  T(\hat{\zeta}^{\tau+1}) u^{\ast} \!+\! W v \leq h(\hat{\zeta}^{\tau+1})\right\} ,
\label{RO12}
\end{equation}
to obtain the adjustable variables $v^{\ast}$ .
The optimal value of the objective function in (\ref{RO12})
is denoted by
$\mathbf{{cost}}_{\mathbf{m},\tau}$. It represents the optimal cost of the problem with the adjustable strategy considered in this section when using method $\mathbf{M_m}$ for determining the nonadjustable variables. If the optimization problem is infeasible we set $\mathbf{{cost}}_{\mathbf{m},\tau}=\infty$.

We  also consider the following  adjustable robust optimization problem according to Theorem \ref{Theorem1}:
\begin{enumerate}
\item[{$\mathbf{M_5}$}]:
\begin{enumerate}
\item Solve the computationally tractable optimization problem (\ref{RO9}) using only the information (\ref{infobis});
\item Solve the optimization problem
\begin{equation}
\min_{\lambda}\varphi^{\tau}(\lambda):=\min_{\lambda}\|\zeta \!- \!\sum_{s=1}^{\tau}\lambda_s(\zeta)\zeta^s\|^2,\ \sum_{s=1}^{\tau}\lambda_s\!=\!1,\ \lambda_s\geq 0,
\label{RO13}
\end{equation}
with  $\zeta=\hat{\zeta}^{\tau+1}$;
\item For each $\tau=\bar{S},\dots,S-1$ define the adjustable variables as in (\ref{v}).
\end{enumerate}
\end{enumerate}
\noindent This method has the advantage that the adjustable variables are obtained by solving only the very simple optimization problem (\ref{RO13}). However, this method only works if the optimal objective function of the latter problem $\varphi^\tau\left(\lambda\left(\hat{\zeta}^{\tau+1}\right)\right)$ is equal to zero. In this case the optimal cost given by this method,
$\mathbf{{cost}}_{\mathbf{5},\tau}$, is equal to the optimal value of the objective function (\ref{RO9}). If $\varphi^\tau\left(\lambda\left(\hat{\zeta}^{\tau+1}\right)\right)>0$  we set $\mathbf{{cost}}_{\mathbf{5},\tau}=\infty$.

Of course, when
$\hat{\zeta}^{\tau+1}$ become available we can also solve the optimization problem
\begin{equation}
\min_{u,v,w}\left\{w:\  p(\zeta)^{T}u + q^{T} v -w \leq - o(\zeta),  T(\zeta) u + W v \leq h(\zeta)\right\}\ ,
\label{RO14}
\end{equation}
 with  $\zeta=\hat{\zeta}^{\tau+1}$. However, this optimization problem, corresponding to the well-known \textit{Wait-and-See} problem (WS) in SP,  determines the optimal values of both the adjustable and nonadjustable variables, while in our setting the  nonadjustable variables have to be determined before the $\hat{\zeta}^{\tau+1}$ become available.
We denote by $\mathbf{{cost}}_{\tau}$ the optimal value of the objective function (\ref{RO14}).

We finally compare methods $\mathbf{M_m},\;\mathbf{m}=\mathbf{1},\ldots,\mathbf{5}$ by considering the aggregated costs:
\begin{equation}
\mathbf{{cost}}_{\mathbf{m}}:= \sum_{\tau=\bar{S}}^{S-1} \mathbf{{cost}}_{\mathbf{m},\tau},\quad \mathbf{m}=\textbf{1},\dots,\textbf{5}\ .
\end{equation}
The same for the wait-and-see problem, WS (\ref{RO14}), $\mathbf{cost}:= \sum_{\tau=\bar{S}}^{S-1} \mathbf{{cost}}_{\tau}\ . $

\section{Supply planning  under uncertainty}
\label{SPUU}
In this section we shortly review the literature  of supply planning under uncertainty and we describe the problem considered to illustrate the efficiency of the scenario based framework methodology.
\subsection{Literature review}
\label{sec:LiteratureReview}

Freight transportation \cite{2} is one of todays' most important logistic activities that influences the performance of many other economic activities. The two most important factors for having high performance level are the economic efficiency and the service quality. The former relates to the facts that those firms that make use of freight transportation service want to move the right amount of  goods at the best cost. The latter highlights the importance of the quality of service, i.e. being able to satisfy clients demand of a good eventually avoiding inventory costs. Bolstering one almost certainly causes the other to suffer: for example, if one lowers his/her inventory to reduce costs,  it will become  difficult to meet varying customer demand. If one increases safety stock to meet peak demands, one could wind up with a great deal of excess inventory with nowhere to sell it, see \cite{14}. Traditionally, two prevailing supply chain strategies have dominated the industry: push and pull. However, new technologies,  have enabled the creation of a third strategy, a hybrid push-pull model that offers the best of both worlds without their corresponding disadvantages, see  \cite{13}.

The problem of transporting goods or resources from a set of supply points (production plants) to a set of demand points (destination factories or customers) is an important component of the planning activity of a manufacturing firm.
A particular case is given by the so-called single-sink transportation problem, in which
a single retailer is served by a set of suppliers. This problem has been extensively studied, in
particular when the total cost is given by the sum of a variable transportation cost and
a fixed charge cost to use the supplier (\cite{AK05,LW97,LSP90,MKB}).
Critical parameters such as customer demands, row material prices and resource capacity are quite uncertain in real problems.
An important issue is then represented by the decision on quantities to acquire and store at each destination factory before actual demands reveal themselves.
This is involved in the tactical planning \cite{2} of the  firm supply chain operations.
In this paper we are dealing with a tactical planning problem of a firm supply chain operations,  where decision on quantities to acquire and store at each destination factory before actual demands reveal themselves
have to be chosen. Detailed information on the capacity of  transportation vehicles, capacity of origin plants  and destination warehouse are given.
See \cite{2} and \cite{G93} for differences among  planning levels, strategic, tactical and operational ones.
Example of strategic models  are those for  designing physical networks and their evolutions,  for location of the main facilities of a plant, for resource acquisition, for defining tariff policies. They usually involve long-term investment decision and strategic policies and they deal with planning at international, national and regional levels.
Tactical planning problems modeling works over a medium term horizon for a rational and efficient allocation of existing resources in order to enhance the performance of the whole system.
They usually incorporate seasonal changes in the data and not daily information: vehicle routing models belong to this class of models.
On the other hand, the notion of cost is central to the operational planning model and it may have different components related to different events that may happen at the decision time (scarcity of goods at the origin, full capacity at destinations, delay, risk, etc.) see \cite{G93}.

The significance of uncertainty has prompted a number of works addressing random parameters in tactical level supply chain planning involving distribution of raw material and products (see for example  \cite{CSS2014,CP96,CL77,CL97,PoTo2003,15,SCL2011,LV02,YL00}). Nowadays firms are engaged in a continuous procurement process and collaboration with their supply chain partners.  The firm regularly orders products from suppliers in a given area according to inventory, forecasted demand and inventory policy \cite{2a}. Firms may negotiate directly with carriers but very often they deal with a third-party logistic service provider (3PL) \cite{10}.
In the recent past, 3PL, also referred to as logistics outsourcing (e.g. \cite{3,9,10,12}), has received considerable attention  due to a steadily increasing number of companies across industry sectors using third-party providers for the management of their logistics operations (e.g. \cite{4,5,6,8}). The functions performed by the third party can encompass the entire logistics process or selected activities within that process \cite{1,7}.
The 3PL combines goods into containers or trucks ensuring shipping usually with a fixed schedule.
The firm needs to book in advance sufficient capacity to satisfy the demand leading to a negotiation with 3PL to reserve the necessary number of trucks.
Usually the agreement between the firm and 3PL is  made \textit{under uncertainty} without full information about the demand of goods.
Extra vehicles must be then purchased at a much higher cost than the initial one.


\subsection{Problem description}
\label{SOMSTP}

We consider a supply planning problem  to optimize vehicle-renting and transportation procurement activities to satisfy demand in several destinations out of several origins.
Uncertainty on the demands and on the costs of extra vehicles is considered.

The logistic system is organized as follows: a set of  suppliers, each of them composed by a set of plants  (origins)  has to satisfy the demand  of a certain good from a set  of factories (destinations)  belonging to the same  producer.
The demand  at  the destination factory  is considered stochastic.
We assume a uniform fleet of vehicles with fixed capacity and we allow only full-load shipments.
Shipments are performed by booking a number  of capacitated vehicles  in advance, before the demand is revealed.
When the demand  becomes known, the number of used vehicles  is determined and there is an option to discount  vehicles booked but not actually used.
The cancellation fee is given as a proportion  of the transportation costs. If the quantity shipped from the suppliers using the booked vehicles is not enough to satisfy the demand at destination factory, residual product  is purchased from an external company  at a higher price, which is considered stochastic as well.
The problem is to determine  the number of vehicles  to book from each plant   of the  set of suppliers  to replenish the good at factory   in order to minimize the total cost for the producer of good, given by the sum of the transportation costs   from origins  to destinations.

\subsection{Notation}
We adopt the following notation.

\label{Not}
\noindent Sets:
\begin{eqnarray*}
\mathscr{K} = \{k:k=1,\ldots ,K\}&,& \textrm{set of suppliers};\\
\mathscr{O}_k = \{i:i=1,\ldots ,O_k\}&,& \textrm{set of plant locations of supplier } k\in \mathscr{K};\\
\mathscr{D} = \{j:j=1,\ldots ,D\}&,& \textrm{set of destination factories};\\
\end{eqnarray*}
\noindent Deterministic parameters:
\begin{eqnarray*}
t_{ijk}&,&  \textrm{unit transportation cost from supplier } i\in\mathscr{O}_k, k\in\mathscr{K} \textrm{ to plant } j\in\mathscr{D}; \\
q &,&  \textrm{vehicle capacity};\\
g_j &,&  \textrm{maximum capacity which can be  booked at the customer } j\in\mathscr{D};\\
v_k &,& \textrm{maximum requirement capacity of supplier }  k\in\mathscr{K};\\
r_k &,& \textrm{minimum requirement capacity of supplier }  k\in\mathscr{K};\\
l_j^{0} &,&  \textrm{initial inventory of product at customer } j\in\mathscr{D}. \\
\alpha &,&  \textrm{discount};\\
\end{eqnarray*}
\noindent Uncertain or Stochastic Parameters:
\begin{eqnarray*}
d_{j} &,&  \textrm{ demand of customer \textit{j}}\in\mathscr{D};\;\\
b_{j}&,&  \textrm{ buying cost from external sources for customer \textit{j}}\in\mathscr{D}. \\
\end{eqnarray*}
\noindent Variables:
\begin{eqnarray*}
x_{ijk}&,& \textrm{number  of  vehicles booked from supplier } i\in\mathscr{O}_k, \  k\in\mathscr{K} \textrm{ to plant } j\in \mathscr{D},\\
&& \textrm{(nonadjustable or first stage decision variables);}\;\\
z_{ijk} &,&  \textrm{number of vehicles actually used from supplier } i\in\mathscr{O}_k,\  k\in\mathscr{K} \textrm{ to plant } j\in \mathscr{D},\\
&&\textrm{ (adjustable or second stage decision variables)};\; \\
y_j &,& \textrm{volume of product to purchase from an external source normalized by $q$, for}\\
&& \textrm{ plant }  j\in \mathscr{D},  \textrm{ (adjustable or second stage decision variables).}
\\
\end{eqnarray*}
To simplify the notation it is convenient to introduce the vectors
\begin{eqnarray*}
&&d=\left(d_j\right)_{j\in\mathscr{D}}\in\mathbb{R}^D,\; b=\left(b_j\right)_{j\in\mathscr{D}}\in\mathbb{R}^D,\; y=\left(y_j\right)_{j\in\mathscr{D}}\in\mathbb{R}^D,\\
&& x={\rm vec}\left(x_{ijk}, i\in\mathscr{O}_k, k\in\mathscr{K},j\in\mathscr{D}\right),\;
z ={\rm vec} \left(z_{ijk}, i\in\mathscr{O}_k, k\in\mathscr{K},j\in\mathscr{D}\right).
\end{eqnarray*}
It is easily seen that $x,z\in\mathbb{R}^m$,  with $\;m=D\sum_{k=1}^KO_k$.
\subsection{Objective function and constraints}
The objective function of our optimization problem is of the form
\begin{equation}\label{obc}
f(x,y,z;b)=f_1(x)+f_2(x,y,z;b),
\end{equation}
where
\begin{equation}\label{obc1}
f_1(x)=q \sum_{k=1}^{K}\sum_{i=1}^{O_k}\sum_{j=1}^D t_{ijk} x_{ijk}
\end{equation}
denotes the
booking costs of the vehicles, while
\begin{equation}\label{obc}
f_2(x,y,z;b)=
   \!\!q\sum_{j=1}^D   b_j\,y_j\!  - \! \alpha q \!\sum_{k=1}^{K}\sum_{i=1}^{O_k}\!\sum_{j=1}^D  t_{ijk}\! \left(x_{ijk}\! -\! z_{ijk} \right),
\end{equation}
represents the cost of recourse
actions, consisting of buying good from external sources ($q y_j$) and canceling unwanted
vehicles. We note that $f_1$ involves only deterministic parameters and nonadjustable, or first stage decision variables while, while $f_2$ also involves uncertain or stochastic  parameters as well as adjustable or second stage decision variables.

We have the following constraints for our optimization problem. Constraint
\begin{equation}\label{c1}
\mathscr{C}_1(x): \qquad  q \sum_{k=1}^{K}\sum_{i=1}^{O_k} x_{ijk} \leq  g_j\ ,\quad  j \in \mathscr{D},
\end{equation}
 guarantees that for each destination $j\in\mathscr{D}$ the number of booked vehicles does not exceed $g_j/q$ inducing thus an upper bound on the total number of
vehicles. We impose the constraint
\begin{equation}\label{c2}
\mathscr{C}_2(z): \qquad r_k \leq q \sum_{i\in O_k}\sum_{j=1}^D z_{ijk} \leq v_k \ ,\quad k \in \mathscr{K},
\end{equation}
to ensures that the number of vehicles serving supplier $k$
does not exceed the production capacity $v_{k}$ of supplier $k\in\mathscr{K}$ and satisfies the lowest requirement capacity $r_k$ established in the contract.
In order to  ensure that the
 $j$-customer's demand is satisfied we also impose the following constraint
\begin{equation}\label{c3}
\mathscr{C}_3(y,z;d): \qquad  l_j^{0} + q \left(\sum_{k=1}^{K}\sum_{i=1}^{O_k} z_{ijk} + y_j \right)- d_j \geq  0\ , \quad  j \in \mathscr{D}.
\end{equation}

The inequality constraint
\begin{equation}\label{zltx}
z\le x,
\end{equation}
guarantees that the number of vehicles actually used
is at most equal to the number booked in advance. We will always impose  nonnegativity constraints on the vectors $y,x,z$,
\begin{equation}\label{c4}
y\ge 0,\; x\ge 0,\;z\ge 0,
\end{equation}
and sometimes integrality constraints on the vectors $x$ and $z$, i.e.,
\begin{equation}\label{c5}
x \in  \mathbb{N}^m,\;\;z \in  \mathbb{N}^m,
\end{equation}
where $\mathbb{N}$ is the set of  natural numbers.

If the vectors $b$ and $d$ were known the optimization problem presented above would be a simple linear programming problem. However, in our application $b$ and $d$ are uncertain. In the next section we will deal with this problem
by using a stochastic programming model and several robust optimization models.
\section{Models for the supply planning problem under uncertainty}
\label{Models}
\subsection{A two stage stochastic optimization model}
\label{SP}
In this section we introduce a two-stage stochastic optimization model for the problem described above. Our approach is based on a set of scenarios
\[
\mathscr{S} = \{s:s=1,\ldots ,S\},
\]
and the following parameters and variables related to it:
\begin{eqnarray*}
P^s &,&  \textrm{probability of scenario } s\in\mathscr{S}; \\
d_{j}^s &,&  \textrm{ demand of customer \textit{j} at scenario } s\in\mathscr{S};\\
b_{j}^s &,&  \textrm{ buying cost from external sources for customer \textit{j} at scenario } s\in\mathscr{S};\\
z_{ijk}^s &,&  \textrm{number of vehicles actually used from supplier } i\in\mathscr{O}_k,\  k\in\mathscr{K} \textrm{ to plant } j\in \mathscr{D},\\
&&\textrm{at scenario } s\in\mathscr{S} \textrm{ (second stage decision variables)};\\
y_j^s&,& \textrm{volume of product to purchase from an external source normalized by } q,\\
&& \textrm{ for plant }  j\in \mathscr{D}, \textrm{at scenario } s\in\mathscr{S} \textrm{ (second stage decision variables)}.
\end{eqnarray*}
\noindent As in subsection~\ref{Not} we introduce the vectors
\begin{eqnarray*}
&&d^s=\left(d_j^s\right)_{j\in\mathscr{D}},\; b^s=\left(b_j^s\right)_{j\in\mathscr{D}},\; y^s=\left(y_j\right)_{j\in\mathscr{D}},\;z^s = \left(z_{ijk}^s\right)_{ i\in\mathscr{O}_k, k\in\mathscr{K},j\in\mathscr{D}}.
\end{eqnarray*}
In the two-stage (one-period) case, we get the following mixed-integer stochastic programming model with recourse:
\begin{eqnarray}\label{objectivefmultistage2}
\mbox{\bf SP}\hspace{1cm}	& &\min_{x,\{y^s,z^s\}_{s\in\mathscr{S}}}\; f_1(x)+  \sum_{s=1}^{S}P^s\left(f_2(x,y^s,z^s,b^s)\right)\\
\label{it1}
&s.t. \quad &\mathscr{C}_1(x),\;\;  \mathscr{C}_2(z^s),\;\mathscr{C}_3(y^s,z^s;d^s),\; s\in\mathscr{S},\\
&&x\in\mathbb{N}^m,\;\;z^s\in\mathbb{N}^m,\; z^s\le x,\;y^s\ge 0,\; s\in\mathscr{S}.\label{it7}
\end{eqnarray}
\noindent  From now on we will refer to problem (\ref{objectivefmultistage2})-(\ref{it7}) as the two-stage stochastic programming problem (SP) for the supply planning problem.

\subsection{Robust optimization models }\label{s:rob}
In this section we introduce several robust formulations (RO) for the problem described above.
These formulations are of particular interest in the case  it is impossible, or not practical, to give reasonable estimates of probability distributions for the random parameters given by the demand  of a certain good at factories and the cost of buying from external sources. Moreover, some LP-relaxations  of RO formulations can be solved in polynomial time and have theoretical guarantees for the quality of the solution none of which is true for the SP formulations.

We consider different selections of the uncertainty set for the objective function and constraints involving the uncertain demands and buying costs. More precisely we assume that $b\in\mathscr{U}_b$ and $d\in\mathscr{U}_d$ for some uncertainty sets $\mathscr{U}_b\subset\mathbb{R}^{D}$ and $\mathscr{U}_d\subset\mathbb{R}^{D}$. For any such uncertainty sets the robust optimization formulation of our problem becomes
}
\begin{eqnarray}
\mbox{\bf RO}\hspace{1cm}		&& \min_{w,x,y,z}\; w
\label{objectiveROg}\\
&\textrm{ s.t. }\; & w \geq f(x,y,z;b),\;\forall b\in\mathscr{U}_{b} \label{itro1g} \\
&&\mathscr{C}_1(x), \mathscr{C}_2(z), \mathscr{C}_3(y,z;d), \;\forall d\in\mathscr{U}_{d}\\
&&x\in\mathbb{N}^m,\;\;z\in\mathbb{N}^m,\; z\le x,\;y\ge 0,
 \label{itro8g}		
\end{eqnarray}
where the auxiliary variable $w$ has been introduced.

\subsubsection{Box uncertainty}
\label{sec:ext-box}
We assume that the cost vector $b$  belongs to the uncertainty set
\begin{equation}
\mathscr{U}_{b, box,L} = \left\{ \bar{b}+\sum_{\ell=1}^L\zeta_\ell b^\ell\;:\; \zeta=[\zeta_1;...;\zeta_L]\in\mathbb{R}^L,\, \|\zeta\|_{\infty}\le 1\right\},
\label{box2e}
\end{equation}
where $b^1,\ldots,b^L$ are vectors representing possible perturbations of the  average  vector cost $\bar{b}$, which is considered fixed.
If we choose $L=D$ and the perturbation vectors
\begin{equation}
\label{choice}
b^\ell=\rho_2F_\ell e^\ell,\; \ell=1,\ldots,D\,,
\end{equation}
where $e^\ell$ is the $\ell-$th vector from the standard basis of $\mathbb{R}^D$
then $y^Tb^\ell=\rho_2F_\ell y_\ell$, with the positive number $F_\ell$ representing the uncertainty scale and $\rho_2>0$  the uncertainty level.

With the choice (\ref{choice}) the uncertainty set (\ref{box2e}) coincides with the simple box
\begin{equation}
\mathscr{U}_{b, box} = \left\{b\in\mathbb{R}^{D}: |b_j - \bar{b}_j|\leq \rho_2 F_{j}, \ j\in\mathscr{D}\right\}.
\label{box2}
\end{equation}
Of course, for other choices of the perturbation vectors we get  different results.

Similarly, we assume that the demand vector $d$ belongs to an uncertainty set of the form
\begin{equation}
\mathscr{U}_{d, box,L} = \left\{ \bar{d}+\sum_{\ell=1}^L\zeta_\ell\, d^\ell\;:\; \zeta=[\zeta_1;...;\zeta_L]\in\mathbb{R}^L,\, \|\zeta\|_{\infty}\le 1\right\},
\label{box1e}
\end{equation}
for given perturbation vectors $d^1,\ldots,d^L$.
As above, it is easily seen that with the choice
\begin{equation}
L=D,\quad d^\ell=\rho_1G_\ell e^\ell,\; \ell=1,\ldots,D\,,
\label{choiced}
\end{equation}
the uncertainty set $\mathscr{U}_{d, box,L}$ reduces to the simple box
\begin{equation}
\mathscr{U}_{d, box} = \left\{d\in\mathbb{R}^{D}: |d_j - \bar{d}_j|\leq \rho_1 G_{j}, \ j\in\mathscr{D}\right\}.
\label{box1}
\end{equation}
We clearly have $\max_{d\in\mathscr{U}_{d, box}}d_j=\bar{d}_j+\rho_1 G_{j}$.
Using the uncertainty sets (\ref{box2}) and $(\ref{box1})$ for the uncertain vectors $b$ and $d$, and  considering the vector $G=\left(G_j\right)_{j\in\mathscr{D}}$, the robust formulation (\ref{objectiveROg})-(\ref{itro8g}) of our problem can be written as the following linear mixed-integer problem:

\begin{eqnarray}
\mbox{\bf RO-box}\hspace{.5cm}	&& \min_{w,x,y,z}\; w
\label{objectiveRObox}\\
&\textrm{ s.t. }\; & w - f(x,y,z;\bar{b})\ge\sum_{j=1}^D  q \rho_2 F_j y_j    \label{itrobox} \\
&&\mathscr{C}_1(x), \mathscr{C}_2(z),  \mathscr{C}_3(y,z;\bar{d}+\rho_1G),\\
&&x\in\mathbb{N}^m,\;\;z\in\mathbb{N}^m,\; z\le x,\;y\ge 0,
 \label{itro8box}		
\end{eqnarray}

We note  that the above model is very conservative. Indeed let us assume that
\begin{equation}\label{zzz}
\zeta_1,\ldots,\zeta_L:\mbox{\rm zero mean independent random variables in the interval } [-1,1] ,
\end{equation}
and let us consider the random vectors
\begin{equation}\label{bdz}
b=b(\zeta)=\bar{b}+\sum_{\ell=1}^D\zeta_\ell\, b^\ell,\quad d=d(\zeta)=\bar{d}+\sum_{\ell=1}^D\zeta_\ell\, d^\ell\, ,
\end{equation}
where the perturbation vectors $b^\ell$ and $ d^\ell$ are given by (\ref{choice}) and (\ref{choiced}) respectively.
Consider a feasible solution of the RO problem (\ref{objectiveRObox})-(\ref{itro8box}). Then
\begin{equation}\label{prob1}
{\mbox{\rm Prob}}_{\zeta\sim P}\left\{\zeta: \,w   \geq f(x,y,z,b(\zeta))\right\}
\!=\!1,
\end{equation}
and
\begin{equation}\label{prob2}
{\mbox{\rm Prob}}_{\zeta\sim P}\left\{\zeta\, : \,\mathscr{C}_3(y,z,d(\zeta))\right\}
=1 ,
\end{equation}
for {\it any} probability distribution $P$ that is compatible with (\ref{zzz}). This certitude of constraints satisfaction will result in a high cost for the optimal solution of the RO problem (\ref{objectiveRObox})-(\ref{itro8box}).

\subsubsection{Box-Ellipsoidal Uncertainty}
\label{sec:BoxEllipsUncertainty}
In this subsection we study the case were the uncertainty set for the buying costs is given by
\begin{equation}
\mathscr{U}_{b,ell} = \left\{ \bar{b}+\sum_{\ell=1}^L\zeta_\ell\, b^\ell\;:\; \zeta=[\zeta_1;...;\zeta_L]\in\mathbb{R}^L,\, \|\zeta\|_{2}\le \Omega\right\},
\label{box2ell}
\end{equation}
and the uncertainty set for the demand vector $d$ is the box (\ref{box1}).

Using the Cauchy-Schwarz inequality, we obtain, for a given $y$
\begin{equation}
\max_{b\in\mathscr{U}_{b,ell}} y^T b=y^T\bar{b}+\max_{\|\zeta\|_2\le \Omega}\sum_{\ell=1}^L\zeta_\ell (y^T b^\ell)=y^T\bar{b}+ \Omega\sqrt{\sum_{l=1}^L (y^T b^\ell)^2}\quad  \ .
\label{b12}
\end{equation}
By choosing the perturbation vectors as in (\ref{choice}) we obtain the following RO model with box and ellipsoidal uncertainty set:
\begin{eqnarray}
\mbox{\bf RO-ell}\hspace{.5cm}	&& \min_{w,x,y,z}\; w
\label{objectiveROell}\\
&\textrm{ s.t. }\; & w - f(x,y,z;\bar{b})\ge\Omega \cdot \sqrt{\sum_{j=1}^D  (q \rho_2 F_j y_j)^2}  \label{itro1ell} \\
&&\mathscr{C}_1(x), \mathscr{C}_2(z),  \mathscr{C}_3(y,z;\bar{d}+\rho_1G),\\
&&x\in\mathbb{N}^m,\;\;z\in\mathbb{N}^m,\; z\le x,\;y\ge 0.
 \label{itro8ell}		
\end{eqnarray}
The nonlinear constraint (\ref{itro1ell}) is a second order cone constraint, so that the above optimization problems is a SOCP. Notice that if we relax the integrality constraints on  $x$ and $z$, the SOCP problem can be solved in polynomial time.

 Consider again the random vectors $b(\zeta)$ and $d(\zeta)$ from (\ref{bdz}). By virtue of Proposition~\ref{prop} it follows that for any feasible solution of the RO problem (\ref{objectiveROell})-(\ref{itro8ell}) we have
\begin{eqnarray}\label{prob1e}
{\mbox{\rm Prob}}_{\zeta\sim P}\left\{\zeta\;  : \;  w \ge f(x,y,z;b(\zeta)) \right\} \ge 1-e^{-\frac{\Omega^2}{2}} ,
\end{eqnarray}
for all probability distributions $P$ that are compatible with (\ref{zzz}).
Since we are using the same box uncertainty for the demand,  (\ref{prob2}) is also satisfied.
From the Cauchy-Schwarz inequality we have
\begin{equation}
\sum_{j=1}^D  q \rho_2 F_j y_j\le\sqrt{D}\sqrt{\sum_{j=1}^D  (q \rho_2 F_j y_j)^2}.
\label{formulaD}
\end{equation}
If follows that if $\Omega\ge\sqrt{D}$ then for any feasible solution of the RO problem (\ref{objectiveROell})-(\ref{itro8ell}) we have the stronger probability result (\ref{prob1}).
This certitude of constraint satisfaction will result in a high cost for the optimal solution of the RO problem (\ref{objectiveROell})-(\ref{itro8ell}).

\subsubsection{Adjustable robust optimization}
\label{sec:adjust}
In the  robust optimization models considered in \ref{sec:ext-box} and \ref{sec:BoxEllipsUncertainty}, all the variables are treated in the same way, while in the two-stage stochastic programming model the variables $x_{ijk}$ are considered first stage decision variables and the variables $y_j$ and $z_{ijk}$ are considered second stage (recourse) variables.
This means that the variables $x_{ijk}$ are to be determined ``here and now'', before the actual data ``reveals itself''.  On the other hand, once the uncertain data are known the variables  $y_j$, $z_{ijk}$ should be able to adjust themselves by means of some decision rules  $Y_j(\cdot)$ and $Z_{ijk}(\cdot)$. Therefore the decision variables $x_{ijk}$ are called nonadjustable  variables while the decision variables $y_j$ and $z_{ijk}$  are called adjustable  variables. In the following we relax the integrality in all the decision variables and we assume that decision rules  $Y_j(\cdot)$ and $Z_{ijk}(\cdot)$ are affine function of their arguments.
In the case with integer adjustable variables  more subtle approaches are needed (e.g. \cite{BG2013,BG2014,HKW2015}).

In developing an adjustable robust optimization model for our problem we will follow the simple model described in \cite{BTEGN09}, where, however, it is assumed that all the coefficients of the adjustable variables are certain. This is not the case for our problem where the coefficients $b_j$ of the adjustable variables $y_j$ are uncertain. We will circumvent this difficulty by assuming as before that the cost vector $b$ belongs to the ellipsoidal uncertainty set $\mathscr{U}_{b,ell}$ (\ref{box2ell}). We have seen that in this case, with the choice (\ref{choice}), our cost constraint can be written under the form (\ref{itro1ell}). This is no longer a linear constraint, but a second order cone constraint. On the other hand we assume that the demand vector $d$ belongs to the scenario-generated uncertainty set
\begin{equation}\label{sce}
\mathscr{U}_{d,\hat{\Delta} } = \left\{\sum_{s=1}^S\lambda_s\hat{d}^s\; : \; \lambda\in\mathscr{L}\right\},
\end{equation}
where,
\begin{equation}\label{sceL}
\mathscr{L}=\left\{\lambda=[\lambda_1;\ldots;\lambda_S]\in\mathbb{R}^S\; : \; \lambda_1\ge 0, \ldots ,\lambda_S\ge 0,\;\sum_{s=1}^S\lambda_s=1\right\},
\end{equation}
and $\hat{\Delta}$ is a given set of scenarios
\begin{equation}\label{sceS}
\hat{\Delta}=\left\{\hat{d}^1,\hat{d}^2,\ldots,\hat{d}^S\right\}.
\end{equation}
In our applications $\hat{\Delta}$ is obtained from historical data.\\
Let us denote by $u=x$ the vector composed of all ``here and now'' decision variables $x_{ijk}$,
and by $v=[y;z]$ the vector composed of all adjustable decision variables $y_j$, $z_{ijk}$.
We consider also the vector of decision rules
\[
V(\cdot)=\left[Y_1(\cdot);\ldots;Y_D(\cdot);\mbox{\rm vec}\left(Z_{ijk}(\cdot),  i \in \mathscr{O}_k,\, k \in \mathscr{K},\,  j \in \mathscr{D}\right)\,\right].
\]
Since decision rules  $Y_j(\cdot)$ and $Z_{ijk}(\cdot)$ are assumed to be affine, so is $V(\cdot)$.
The deterministic constraints of our problem are:
\begin{equation}\label{dete}
\mathscr{C}(u,v):\quad\mathscr{C}_1(x),\; \mathscr{C}_2(z),\;z\le x,\;x\ge0,\;y\ge 0,\;z\ge 0,
\end{equation}
while our uncertain constraints are given by:
\begin{equation}\label{uns}
\widetilde{\mathscr{C}}(u,v,w;b,d):\quad w \ge f(x,y,z;b),\; \mathscr{C}_3(y,z;d).
\end{equation}
With the above notation our uncertain problem can be written as:
\[
\mathscr{R}=\min_{w,u,v}\left\{w\; :\; \mathscr{C}(u,v),
\widetilde{\mathscr{C}}(u,v,w;b,d),
\forall  b\in\mathscr{U}_{b,ell} ,\forall  d\in\mathscr{U}_{d,\hat{\Delta} }\right\}.
\]
For the choice (\ref{choice}) the infinite set of constraints $ w \ge f(x,y,z;b), \forall  b\in\mathscr{U}_{b,ell}$ reduces to the deterministic single second-order cone constraint (\ref{itro1ell}). In this case the uncertain constraints of our problem become
\begin{equation}\label{unsh}
\widehat{\mathscr{C}}(u,v,w;d):\quad w - f(x,y,z;\bar{b})\ge
		\Omega\sqrt{\sum_{j=1}^D  (q \rho_2 F_j y_j)^2 },\;\; \mathscr{C}_3(y,z;d).
\end{equation}
Therefore, our uncertain problem can be written equivalently, (see (\ref{RO8})) as:
\begin{equation}\label{rr}
\mathscr{R}\quad :\quad \min_{w,u,v}\left\{w\; :\; \mathscr{C}(u,v),\widehat{\mathscr{C}}(u,v,w;d),\;
\forall  d\in\mathscr{U}_{d,\hat{\Delta} }\right\}.
\end{equation}
According to (\ref{RO7}) the adjustable version of this problem is
\begin{equation}\label{aa}
\mathscr{A}\quad : \quad \min_{w,u,V(\cdot)}\left\{w\; :\; \mathscr{C}(u,V(d)),
\widehat{\mathscr{C}}(u,V(d),w;d)
,\forall  d\in\mathscr{U}_{d,\hat{\Delta} }\right\},
\end{equation}
where the minimum is taken for all decision rules $V(\cdot)$ that are affine functions of their arguments.
According to Theorem \ref{Theorem1} we  show below that this adjustable version is equivalent to the following tractable optimization problem (see (\ref{RO9}))
\begin{equation}\label{qq}
\mathscr{Q}\quad : \quad \min_{w,u,\left\{v^s
\right\}_1^S}\left\{w\; :\; \mathscr{C}(u,v^s),
\widehat{\mathscr{C}}(u,v^s,w;\hat{d}^s),\;
 s=1,\ldots,S \right\}.
\end{equation}
The equivalence is understood in the sense that the optimal values of $\mathscr{A}$ and $\mathscr{Q}$ are
equal and that any feasible solution of $\mathscr{Q}$ can be augmented to a feasible solution of $\mathscr{A}$.
More specifically let $\hat{w},\hat{u},\left\{\hat{v}^s
\right\}_1^S$ be a feasible solution of $\mathscr{Q}$ and consider a vector $d\in\mathscr{U}_{d,\hat{\Delta} }$. Then there is a vector $\lambda(d)=[\lambda_1(d);\ldots;\lambda_S(d)]\in\mathscr{L}$ such that
\begin{equation}\label{dd}
d=\sum_{s=1}^S\lambda_s(d)\hat{d}^s,
\end{equation}
and the adjustable variables are defined by the decision rule
\begin{equation}\label{vv}
v=\hat{V}(d):=\sum_{s=1}^S\lambda_s(d)\hat{v}^s.
\end{equation}
The feasibility of $\hat{w},\hat{u},\left\{\hat{v}^s
\right\}_1^S$  means that the following constraints are satisfied
\begin{equation}\label{feas}
 \mathscr{C}(\hat{u},\hat{v}^s),
\widehat{\mathscr{C}}(\hat{u},\hat{v}^s,\hat{w};\hat{d}^s),\;
 s=1,\ldots,S.
\end{equation}
From (\ref{c1}),  (\ref{c2}) and  (\ref{dete}) it follows that the constraints $\mathscr{C}(u,v)$  are linear inquality constraints, while $\widehat{\mathscr{C}}(u,v,w;d)$ are convex constraints,  because the left hand side of the inequality in (\ref{unsh}) is linear while the right-hand-side can be written as $\| H y\|$ where $H$ is the diagonal matrix with entries $\Omega q \rho_2 F_j$. Therefore if (\ref{feas}) is satisfied, then
$\hat{w},\hat{u},\hat{V}(\cdot)$ is feasible for $\mathscr{A}$. In particular, this proves that the optimal solution of $\mathscr{A}$ is less than or equal to the optimal solution of $\mathscr{Q}$. To prove the reverse inequality, we remark  that
if $(w,u,V(\cdot))$ is feasible for $\mathscr{A}$, then $w,u,\left\{V(\hat{d}^s)
\right\}_1^S$ is clearly feasible for $\mathscr{Q}$.

In conclusion, in order to solve our adjustable robust optimization model we first find the optimal solution
\begin{equation}
x^*=\left(x^*_{ijk}\right),\;y^{*\,s}=\left( y^{*\,s}_j\right),\;z^{*\,s}=
\left(z^{*\,s}_{ijk}\right), \quad	s \in \mathscr{S}
\end{equation}
of the tractable second order cone optimization problem
\begin{eqnarray}
\label{tt}
\mbox{\bf trSOCP} &&\min_{x,\{y^s,z^s\}_{s\in\mathscr{S}}}\;
\quad w \\
&s.t. \quad & w - f(x,y^s,z^s;\bar{b})\ge\Omega \cdot \sqrt{\sum_{j=1}^D  (q \rho_2 F_j y_j^s)^2},\; 	s \in \mathscr{S},
\label{it0a}\\
&&\mathscr{C}_1(x),\;\;  \mathscr{C}_2(z^s),\;\mathscr{C}_3(y^s,z^s;\hat{d}^s),\; s\in\mathscr{S},\\
&&0\le z^s\le x,\;y^s\ge 0,\; s\in\mathscr{S}.\label{it7a}
\end{eqnarray}
When the uncertain demand $d$ reveals itself we try to find a vector \\ $\lambda(d)=[\lambda_1(d);\ldots;\lambda_S(d)]\in\mathscr{L}$ satisfying (\ref{dd}) by solving the following optimization problem in $\lambda$
\begin{equation}\label{lala}
\min_{\lambda\in\mathscr{L}}\;\varphi(\lambda):=\min_{\lambda\in\mathscr{L}}\sum_{j=1}^D\left(d_j-\sum_{s=1}^S\lambda_s\hat{d}^s_j\right)^2\, .
\end{equation}
The  optimal value of the objective function is equal to zero, i.e., $\varphi(\lambda(d))=0$, if and only if $d\in\mathscr{U}_{d,\hat{\Delta} }$. In this case the adjustable variables are given by
\begin{equation}\label{adjv}
y_j^*=\sum_{s=1}^S\lambda_s(d)y^{*\,s}_j, \;
z^{*}_{ijk}=\sum_{s=1}^S\lambda_s(d) z^{*\,s}_{ijk} \quad i \in \mathscr{O}_k,\, k \in \mathscr{K},\,  j \in \mathscr{D}.
\end{equation}
From the above considerations it follows that
\begin{equation}\label{solst}
x^*_{ijk},\; y^{*}_j,\;z^{*}_{ijk}, \quad i \in \mathscr{O}_k,\, k \in \mathscr{K},\,  j \in \mathscr{D},
\end{equation}
is an optimal solution of the robust adjustable optimization problem (\ref{aa}). This is no longer the case  when
$\varphi(\lambda(d))>0$. As noted in Section \ref{sec:AScenarioBasedFrameworkForComparisonOfModellingApproachesUnderUncertainty}, the scenario-generated uncertainty set $\mathscr{U}_{d,\hat{\Delta} }$ from (\ref{sce}) is usually ``too small'' to be of much
interest. Nevertheless, as shown in the next section the values
\begin{equation}\label{solstx}
x^*_{ijk},\quad i \in \mathscr{O}_k,\, k \in \mathscr{K},\,  j \in \mathscr{D},
\end{equation}
could be used to find a solution of our supply planning problem even when $d\notin\mathscr{U}_{d,\hat{\Delta} }$.

\subsubsection{Scenario based framework for the supply planning problem under uncertainty}
\label{sec:adj}

In this subsection we compare the performance of the various methods in the scenario based framework described in Section \ref{sec:AScenarioBasedFrameworkForComparisonOfModellingApproachesUnderUncertainty}. For our specific application those methods become:
\begin{enumerate}
\item[{$\mathbf{M_1}$}] is the  stochastic optimization problem SP (\ref{objectivefmultistage2})--(\ref{it7});
\item[{$\mathbf{M_2}$}] is the RO problem with box constraints RO-box (\ref{objectiveRObox})--(\ref{itro8box});
\item[{$\mathbf{M_3}$}] is the RO problem with ellipsoidal constraints RO-ell (\ref{objectiveROell})--(\ref{itro8ell});
\item[{$\mathbf{M_4}$}] is the second order cone  optimization problem trSOCP (\ref{tt})--(\ref{it7a}).
\end{enumerate}
As in (\ref{sceS}), let $\hat{\Delta} = \left\{\hat{d}^1,\hat{d}^2,\ldots,\hat{d}^{S}\right\}$
be a given  set of demand scenarios of dimension $D$ from historical data. We consider  a set of indices (\ref{bs}).
For each  $\tau=\bar{S},\dots,S-1$ we compute the quantities
$\bar{d}^{\tau}=\frac{1}{\tau}\sum_{s=1}^{\tau} \hat{d}^s,\quad \rho_1^{\tau} G_j^{\tau} = \max_{s=1,\dots,\tau} |\hat{d}_j^s -\bar{d}_j^{\tau}|, \; j\in\mathscr{D}.$

If we do not have historical data for the cost vectors, but we know a $D$-dimensional vector $\bar{b}$ of average costs, we can generate a set of vectors
$\hat{B}=\left\{\hat{b}^1,\hat{b}^2,\ldots,\hat{b}^{S}\right\},$ with components  $\hat{b}_j^s$ obtained by
sampling from a uniform distribution  in the interval $\left[\bar{b}_j - \sigma \cdot \bar{b}_j ,\bar{b}_j + \sigma \cdot \bar{b}_j \right]$ with a  given deviation level of $\sigma$.
For each $\tau=\bar{S},\dots,S-1$ we  compute  $\rho_2^{\tau} F_j^{\tau} = \max_{s=1,\dots,\tau} |\hat{b}_j^s -\bar{b}_j^{\tau}|, \; j\in\mathscr{D}.$
For each $\tau=\bar{S},\ldots,S-1$ and for each of the methods $\mathbf{M_1},\,\ldots,
\mathbf{M_4}$ described in Section~\ref{sec:AScenarioBasedFrameworkForComparisonOfModellingApproachesUnderUncertainty} we  find the optimal solution of the corresponding optimization problem using only the information contained in the vectors
\begin{equation}\label{infotransp}
\hat{d}^1,\hat{d}^2,\ldots,\hat{d}^{\tau},\quad \hat{b}^1,\hat{b}^2,\ldots,\hat{b}^{\tau}.
\end{equation}
From this optimal solution we obtain the nonadjustable decision variables $x^*$.
They are determined by using only the information contained in (\ref{infotransp}). Assume now that the vectors $\hat{d}^{\tau+1},\hat{b}^{\tau+1}$ become available.
Then we can solve the optimization problem
\begin{eqnarray}\label{o1}
&&	 \min_{y,z}  f\left(x^*,y,z;\hat{b}^{\tau+1}\right)\label{o1}\\
&\textrm{ s.t. }\; &  \mathscr{C}_2(z),  \mathscr{C}_3(y,z;\hat{d}^{\tau+1}),\\
&&z\in\mathbb{N}^m,\; z\le x^*,\;y\ge 0.\label{i7}
\end{eqnarray}
to obtain the adjustable variables $y^{\ast}, z^{\ast}$.

We will also consider method $\mathbf{M_5}$ which now consists in the following steps:
\begin{enumerate}
\item[{$\mathbf{M_5}$}]:
\begin{enumerate}
\item Solve the optimization problem (\ref{tt}) -(\ref{it7a}) using only the information (\ref{infotransp});
\item Solve the optimization problem (\ref{lala}) with $S=\tau$ and $d=\hat{d}^{\tau+1}$;
\item If $\varphi(\lambda(d))=0$, define the adjustable variables as in (\ref{adjv}).
\end{enumerate}
\end{enumerate}

\section{Numerical results}
\label{sec:numericalresults}

In this section we discuss  numerical results for  stochastic and robust modelling approaches applied to the supply planning problem described in Sections \ref{SPUU} and \ref{Models}.
The instance considered comes from a real case of gypsum replenishment in Italy provided by the primary Italian cement producer.
The logistic system is composed of 24 suppliers with multiple plants and 15 destinations.
The models  aim to find, for each plant of the 24 suppliers, the number of vehicles to be booked for replenishing gypsum  in order to minimize the total cost.

Deterministic and stochastic parameters values  for the problem are reported in Tables \ref{tab1}-\ref{tab2}-\ref{tabupper} in the Annex. Table~\ref{tab1} lists the set of suppliers $\mathscr{K}$ and the set of their  plants $\mathscr{O}_k, \ k \in\mathscr{K}$. The list of destinations (cement factories) are shown in Table~\ref{tab2} with relative emergency costs, unloading capacities and  variation of demand $\rho_1 G_j $, and buying costs $\rho_2 F_j $, $j\in\mathscr{D}$.

Table~\ref{tabupper} refers to the minimum and maximum requirement capacity of supplier $k \in\mathscr{K}$.
We suppose to have an initial inventory level $l_j^0=0$  at all the destinations  $j\in\mathscr{D}$ and to use vehicles with fixed capacity  $q=31$ tons.
The cancellation fee $\alpha$ is fixed to the value of 0.7.


Equiprobable scenarios of gypsum demand $\hat{\Delta}  = \left\{\hat{d}^1,\hat{d}^2,\ldots,\hat{d}^{S}\right\}$ are built on historical data, using all the weekly values in March, April, May and June of the years 2011, 2012 and 2013 for a total number of $S=48$ realizations.
Having no historical data for the buying costs, equiprobable scenarios  $\hat{B}=\left\{\hat{b}^1,\hat{b}^2,\ldots,\hat{b}^{S}\right\}$ have been generated sampling from a uniform distribution  in the interval\\ $\left[\bar{b}_j - \sigma \cdot \bar{b}_j ,\bar{b}_j + \sigma \cdot \bar{b}_j \right]$ with a  given deviation level of $\sigma=20\%$. In order to make a fair comparison between RO and SP methodologies, the same deviation level is used also  for the RO methodology.



All the models  (SP and  RO) are modeled in AMPL and solved using the CPLEX $12.5.1.0.$ solver or MOSEK 7.1.0.58 solver for RO model with ellipsoidal constraint. The computations have been performed on a 64-bit machine with 12 GB of RAM and a 2.90 GHz Intel i7 processor. Summary statistics of the adjusted problems derived for our test-case are reported in Table \ref{statistics}.

\begin{table}[ht!]
\tiny{
	\centering
		\begin{tabular}{l|l|l|l}
	
			 &  SP $(S=48)$ & RO-box & RO-ell \\
			\hline
	MIP simplex iterations & 16609	 &  83 & -\\
	number of variables     & 24240 & 976 & 976\\
	number of integer variables     & 23520 & 960 & -\\
	number of linear constraints  & 24818 & 521  & 518  \\
	number of non linear constraints  & - & - & 1  \\
	CPU time (s) &	16.19	&  0.109 & 0.2028
		\end{tabular}
		\caption{Summary statistics of the solution approaches (SP versus static RO).}
		\label{statistics}
		}
\end{table}



\noindent First-stage solutions of the stochastic mixed-integer model (SP) are reported in Table \ref{tab5}. Results show that the demand of all 15 destinations is satisfied by ordering vehicles from the set of suppliers for a total of 141 vehicles and a total cost of $107\,244.67$.
The SP model forces the company to buy from external sources only in high demand scenarios as recourse action.

\begin{figure}[ht!]
\centering
\includegraphics[width=\textwidth]{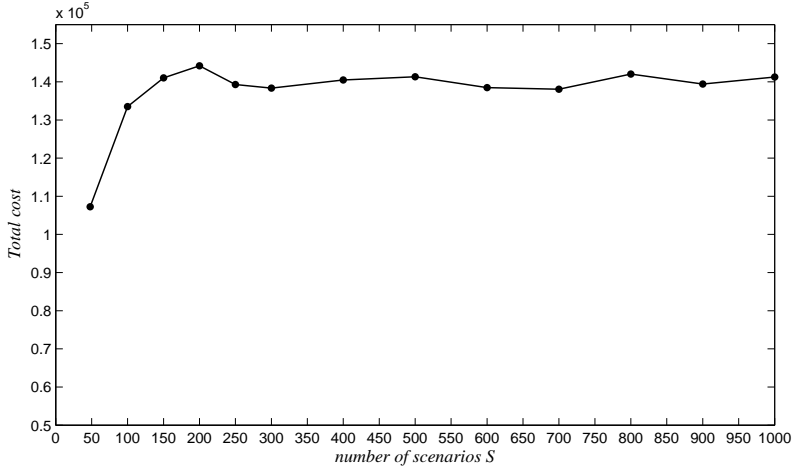}
\vskip3truemm\noindent
\caption{In-sample stability of SP model with all continuous variables.}
\label{convergence}
\end{figure}

Similar results are obtained also for larger number of scenarios. Figure \ref{convergence} shows the \textit{in-sample stability} of  the optimal SP function value  for an increasing number of scenarios \cite{KW}. For  $S\geq 200$  the stability is reached, where the scenarios have been generated by a uniform distribution in the range provided by the minimum and maximum values from historical data.


\begin{figure}[ht!]
\centering
\includegraphics[width=\textwidth]{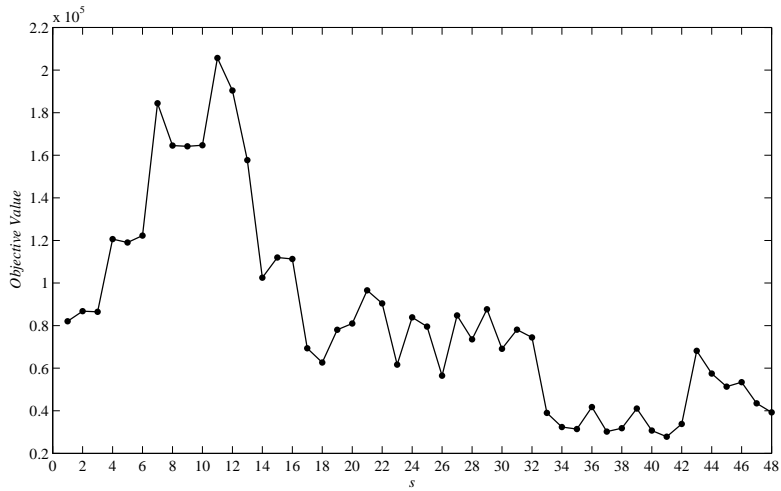}
\vskip3truemm\noindent
\caption{Optimal objective function values of deterministic programs with mixed-integer variables with full knowledge of demand and buying cost. The costs refer to the 48 scenarios obtained from historical data.}
\label{solutions1}
\end{figure}

Before to proceed with the comparison between SP and RO models, we analyze the advantage of having the information about future demand and buying cost. This is provided by the well-known
\textit{Expected Value of Perfect Information} EVPI, \cite{BL11},\cite{MAB2014}:
\begin{equation}
EVPI = SP - WS = 107\,244.67 - 84\,472.21 = 22\,772.46 \
\end{equation}
given by the difference between the stochastic cost and the wait-and-see cost $WS$.
Figure \ref{solutions1} shows the objective function values of the deterministic problems solved separately over each scenario $(\hat{d}^s,\hat{b}^s)\in \hat{\Delta} \times \hat{B}.$

Since the gypsum demand is highly affected by the economic conditions of the public and private medium and large-scale construction sector, a reliable forecast and reasonable estimates of probability distributions are difficult to obtain.
This is the main reason that leads us to consider also  robust optimization approaches.
In the following we consider:
\begin{itemize}
	\item[-] static approaches with uncertainty parameters belonging respectively  to box, ellipsoidal uncertainty sets or mixture of them;
	\item[-] dynamic approaches, via the concept of  tractable adjustable robust counterpart and scenario based framework.
\end{itemize}

We first consider the robust optimization model with box uncertainty, (\ref{objectiveRObox})-(\ref{itro8box}).
Boxes for demand and buying cost are built in compliance with the 48 historical data used to generate SP scenarios. See again Table \ref{tab2} for variation of demand $\rho_1 G_j $  and buying costs $\rho_2 F_j $, $j\in\mathscr{D}$.

The static RO-box approach  is very conservative having a total cost about $365\%$ larger than the SP expected cost. 
Table \ref{tabROybox} (column 2)  reports the optimal solutions of the mixed-integer robust problem with box constraints given by the  normalized volume $y_{j}$ for  destination plant $j\in\mathscr{D}$, while Table  \ref{tabRO}  gives number of booked vehicles $x_{ijk}$ from plant $i\in\mathscr{O}_k$ of supplier $k$ to destination $j\in\mathscr{D}$.

\begin{table}[ht!]
\centering
{\tiny
\caption{Optimal normalized volume $y_{j}$ for  destination plant $j\in\mathscr{D}$  for the static robust mixed-integer optimization model with box constraints (column 2), of the corresponding continuous version (column 3) and  box-ellipsoidal optimization model with $\Omega=2.75$ (column 4).}
 \label{tabROybox}
\begin{tabular}[l]{@{}l|l|l|l}
& mixed-integer RO box  model &  continuous RO box  model & continuous RO-ell  model \\
\textbf{Destination}  &  $y_{j}$ & $y_{j}$ & $y_{j}$\\
\hline
dest1  &   0.643847  & 0 & 4.36  \\
dest2      &  17.2759 & 17.1097 & 9.12 \\
dest3           &  29.9154 & 29.9154 & 12.60\\
dest4          &  8.7287 & 8.5029 & 8.50 \\
dest5            &  0 &   0 & 2.18\\
dest6            &  0.705345 &   0 & 6.94\\
dest7   &  0.595318 &  0 & 5.05\\
dest8      &  0.931387  &   0 & 1.64\\
dest9             & 15.7849  & 13.2951 & 11.13\\
dest10        & 0.542244 & 0 & 0.48\\
dest11              &  1.83631 & 1.80405  & 9.15\\
dest12         & 22.1082 & 22.1082   & 22.10 \\
dest13       &  0  &   0  & 0\\
dest14           &   6.09235 &  5.513  & 9.86\\
dest15       & 23.3377 & 23.3333 & 11.70
\end{tabular}
}
\end{table}

\begin{table}[ht!]
\centering
{\tiny
\caption{Optimal first-stage solutions of the stochastic model SP with mixed-integer variables. The table shows the optimal number of booked vehicles $x_{ijk}$ from plant $i\in\mathscr{O}_k$ of supplier $k$ to destination $j\in\mathscr{D}$.}
 \label{tab5}
\begin{tabular}[l]{@{}llll}
\textbf{Destination} $j\in\mathscr{D}$ & \textbf{Supplier} $k\in\mathscr{K}$ & \textbf{Plant} $i\in\mathscr{O}_k$ & $x_{ijk}$\\
\hline
dest1 & suppl1 &                 plant3   &             9 \\
\hline
dest2    & suppl1               & plant3      &              22 \\
dest2    & suppl3                   & plant8      &           1 \\
dest2    & suppl4         & plant4      &           2 \\
dest2    & suppl6 & plant10 & 4 \\
dest2    & suppl12                   & plant14       &       3 \\
dest2    & suppl16 & plant16 & 3 \\
\hline
dest3         & suppl1  & plant5            &  8 \\
dest3         & suppl3 & plant8  &   2 \\
dest3         & suppl6 & plant10  &   8 \\
dest3         & suppl15 & plant8  &   3 \\
\hline
dest4       & suppl1  & plant5  &    1 \\
dest4       & suppl3  & plant8 &    1 \\
dest4       & suppl13  & plant15 &    3 \\
dest4       & suppl20  & plant5 &    1 \\
\hline
dest5         & suppl20 & plant5   &  1  \\
\hline
dest6         &  suppl10              & plant12          &   2 \\
dest6         & suppl14             & plant12          &   3 \\
dest6         & suppl24              & plant12          &   2 \\
\hline
dest7 & suppl1  & plant5 &    4 \\
dest7 & suppl3  & plant8 &    1 \\
dest7 & suppl20 & plant5 &    1 \\
\hline
dest8   & suppl5                & plant9    &       1 \\
dest8   & suppl9 & plant11    &       3 \\
dest8   & suppl18                & plant9    &       1 \\
dest8   & suppl22                & plant19    &              3 \\
\hline
dest9          & suppl2             & plant7          &   6 \\
dest9          & suppl5             & plant9          &   3 \\
dest9          & suppl9             & plant11          &   1 \\
dest9          & suppl10             & plant12          &   1 \\
dest9          & suppl23              & plant7          &   1 \\
dest9          & suppl24             & plant12          &   1 \\
\hline
dest10      & suppl1            & plant3                  &   6 \\
\hline
dest11      & suppl11 & plant13   &  7 \\
\hline
dest12      & suppl1   & plant5 &    4 \\
dest12      & suppl11  & plant13 &    2 \\
dest12      & suppl17  & plant17 &    3 \\
dest12      & suppl20  & plant5 &    1\\
\hline
dest13    & suppl5  & plant9                &  3 \\
dest13    & suppl9    & plant11                &  3  \\
dest13    & suppl18   & plant9                &  2 \\
\hline
dest14         & suppl1    & plant5 &   3 \\
dest14         & suppl19 & plant5 &   3 \\
dest14         & suppl20 & plant5  &   1 \\
\hline
dest15    & suppl8             & plant2   &   3 \\
dest15    & suppl10             & plant12          &   1 \\
dest15    & suppl21             &  plant18          &   3 \\
dest15    & suppl24             & plant12          &   1 \\
\hline
 \end{tabular}
}
 \end{table}

\begin{table}[ht!]
\centering
{\tiny
\caption{Optimal solution of the RO-box model. The table shows the optimal number of booked vehicles $x_{ijk}$ from plant $i\in\mathscr{O}_k$ of supplier $k$ to destination $j\in\mathscr{D}$.}
 \label{tabRO}
\begin{tabular}[l]{@{}llll}
\textbf{Destination} $j\in\mathscr{D}$ & \textbf{Supplier} $k\in\mathscr{K}$ & \textbf{Plant} $i\in\mathscr{O}_k$ & $x_{ijk}$\\
\hline
dest1 & suppl1               &  plant3   &                 13 \\
\hline
dest2    & suppl1              & plant3      &              29 \\
dest2    & suppl4                   & plant4       &                2 \\
dest2    & suppl6 & plant10 & 12\\
dest2    & suppl12                & plant14      &              3 \\
dest2    & suppl16 & plant16 & 3 \\
\hline
dest4       & suppl13  & plant15 &    3 \\
\hline
dest5         & suppl3 & plant8   &  2  \\
\hline
dest6         & suppl7               & plant1          &   3 \\
dest6         & suppl10               & plant12          &   3 \\
dest6         & suppl14              & plant12          &   3 \\
dest6         & suppl24            & plant12          &   3 \\
\hline
dest7 & suppl3   & plant8 &    4 \\
dest7 & suppl15   & plant8 &    3 \\
dest7 & suppl20   & plant5 &    2 \\
\hline
dest8   & suppl5  & plant9   &              5 \\
dest8   & suppl9               & plant11    &              4 \\
dest8   & suppl22                & plant19    &              3 \\
\hline
dest9          & suppl8             & plant2         &   3 \\
dest9          & suppl21             & plant18         &   1 \\
\hline
dest10      & suppl1              & plant3                  &   6 \\
dest10      & suppl17      &  plant17 & 3 \\
\hline
dest11           & suppl11 & plant13   &  9 \\
\hline
dest13    & suppl2    & plant7                &  6 \\
dest13    & suppl5    & plant9                &  1 \\
dest13    & suppl18    & plant9                &  3 \\
dest13    & suppl23    & plant7                &  1 \\
\hline
dest14         &  suppl19  & plant5 &   3 \\
dest14         &  suppl20  & plant5  &   1
\end{tabular}
}
\end{table}

Due to the  upper bound constraint on the number of booked vehicles  and largest demand to be satisfied, the model forces the company to buy from external sources. This happens for almost all the destinations, with exception of dest5 and dest13 where the demand is fully satisfied by the orders $x_{ijk}$. On the other hand, the demand of dest3, dest12 and dest15  are  satisfied only by external orders $y_j>0$ with a consequent larger cost.

The choice of the box uncertainty set is preferable only if the feasibility of all the constraints is highly required which would be the case if  the company prefers to keep the same contract with the suppliers (or 3PL) and still being immunized against every possible realization of random demand in the box.

One can try to use a different uncertainty set in order to get a less conservative outcome.  We show the results obtained by using a box-ellipsoidal uncertainty set which, as mentioned in Section \ref{sec:BoxEllipsUncertainty}, requires the solution of a second-order cone program with integer constraints.
The current MOSEK solver version 7.1.0.58 finds that the robust mixed-integer second-order cone  formulation RO-ell described in   Section \ref{sec:BoxEllipsUncertainty} is infeasible. Although the problem may be feasible, MOSEK cannot find it.
This is the main reason that made us relax the integrality constraints   on $x_{ijk}$ and $z_{ijk}$ for all  second-order cone  formulations.
To have a fair comparison we have then relaxed the integrality constraints in all the other formulations.

\begin{figure}[ht!]
\centering
\includegraphics[width=\textwidth]{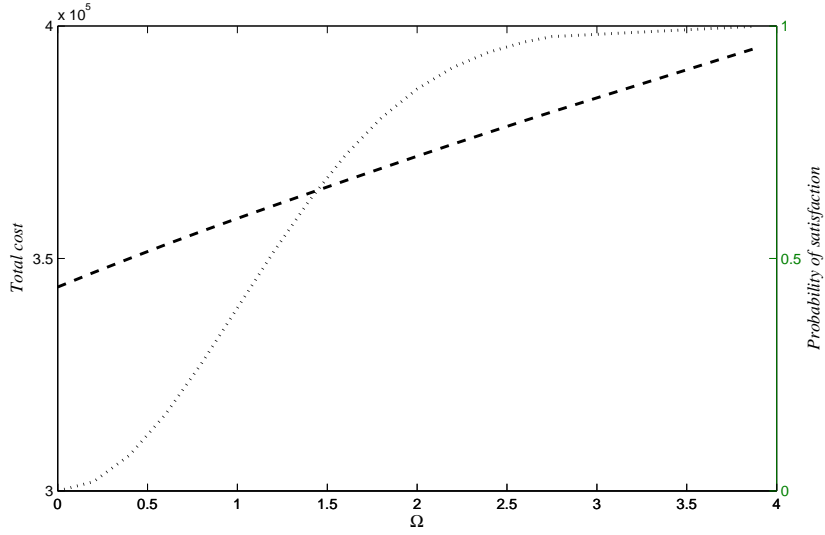}
\vskip3truemm\noindent
\caption{Total cost given by RO-ell with all continuous variables for increasing values of $\Omega$  (dashed line) and probability of satisfaction (dotted line) of constraint (\ref{itro1ell}).}
\label{ellipsoidfigure}
\end{figure}

Total costs of the static RO-ell model are shown in Figure \ref{ellipsoidfigure} (dashed line) with the probability of satisfaction (dotted line) of the second order cone constraint (\ref{itro1ell})   for increasing values of the parameter $0\leq\Omega\leq 3.873$ (see Proposition \ref{prop}).  Notice that the strongest probability satisfaction (\ref{prob1}) is verified  for $\Omega\ge \sqrt{15}=3.873$ since in the test-case considered the number of destinations is $D=15$. The results show that for $\Omega=0$, the total cost of the static RO-ell approach is 343\,849.19, the same than the static RO-box model with average buying cost $\bar{b}_j$ and box constraint requirement for the demand $d_j$.
 As  $\Omega$ increases to 2.75 the total cost reaches approximately the same  value, $381\,520.16$, for  the box model case  with a probability of constraint satisfaction close to one.
Notice that  the optimal cost of the RO-ell model  is only 114.74 lower than the box model, where the probability of constraint satisfaction is exactly  one.
On the other hand, the solution of the static RO-ell model with $\Omega=\sqrt{D}= 3.873$ is guaranteed to satisfy the second order cone constraint with probability one.
This fact has been remarked at the end of Section~\ref{sec:BoxEllipsUncertainty}. For $\Omega=\sqrt{D}= 3.873$ the ellipsoidal uncertainty set $\mathscr{U}_{b,ell}$, given by (\ref{box2ell}) and (\ref{choice}), includes the box uncertainty set $\mathscr{U}_{b, box}$ defined in (\ref{box2}). Therefore the static RO-ell gives a  cost that is $13\,575$ larger than RO-box.
CPU times are negligible (0.2028 CPU seconds) compared to the SP problem with 200 scenarios (114.13 CPU seconds).

For a comparative analysis, Table \ref{tabROell} and Table \ref{tabROybox} (column 4) report the solution variables $x_{ijk}$ and $y_j$ in the case of the static box-ellipsoidal robust approach with $\Omega=2.75$ whereas Table \ref{tabROcont} and Table \ref{tabROybox} (column 3) refer to the optimal solutions of the static box case with continuous variables. While the two approaches have approximately the same total costs, their solution strategies show some  differences: the box-ellipsoidal solution does not make any order except for dest12, deciding to satisfy their maximum demand by external sources $y_{12}=22.10$. On the other side the box solution does not make any order both for dest3  and dest12, and buys from external sources.
The box solution tries to satisfy the demand of dest1, dest5, dest6, dest7, dest8, dest10,  and dest13 only by booking vehicles $x_{ijk}$ from the set of suppliers while the box-ellipsoidal solution requires for all the destinations, with exception of dest13, to buy from external sources.

\begin{table}[ht!]
\centering
{\tiny
\caption{Optimal solution given by RO-box  with all continuous variables. The table shows the optimal number of booked vehicles $x_{ijk}$ from plant $i\in\mathscr{O}_k$ of supplier $k$ to destination $j\in\mathscr{D}$.}
 \label{tabROcont}
\begin{tabular}[l]{@{}llll}
\textbf{Destination} $j\in\mathscr{D}$ & \textbf{Supplier} $k\in\mathscr{K}$ & \textbf{Plant} $i\in\mathscr{O}_k$ & $x_{ijk}$\\
\hline
dest1 & suppl1              &  plant3   &                 13.6438 \\
\hline
dest2    & suppl1              & plant3      &              28.4268 \\
dest2    & suppl4                   & plant4       &                 2.12903 \\
dest2    & suppl6 & plant10 & 12.1587 \\
dest2    & suppl12                & plant14      &              3.22581 \\
dest2    & suppl16 & plant16 & 3.22581 \\
\hline
dest4       & suppl13  & plant15 &    3.22581 \\
\hline
dest5         & suppl3 & plant8   &  1.95436  \\
\hline
dest6         & suppl7               & plant1          &   3.22147 \\
dest6         & suppl10               & plant12          &   3.03226 \\
dest6         & suppl14              & plant12          &   3.22581 \\
dest6         & suppl24            & plant12          &   3.22581 \\
\hline
dest7 & suppl3   & plant8 &    4.49725 \\
dest7 & suppl15   & plant8 &    3.22581  \\
dest7 & suppl19   & plant5 &    1.87226 \\
\hline
dest8   & suppl5  & plant9   &              4.86687 \\
dest8   & suppl9               & plant11    &              4.83871 \\
dest8   & suppl22                & plant19    &              3.22581 \\
\hline
dest9          & suppl2                    & plant7         &    0.425199 \\
dest9          & suppl8             & plant2         &   3.22581 \\
dest9          & suppl21             & plant18         &   3.22581 \\
dest9          & suppl23                         & plant7        &   1.6129 \\
\hline
dest10      & suppl1                 & plant3                  &   6.31644 \\
dest10      & suppl17      &  plant17 & 3.22581 \\
\hline
dest11           & suppl11 & plant13   &  9.03226 \\
\hline
dest13    & suppl2    & plant7                &  6.02641 \\
dest13    & suppl5    & plant9                &  1.58474 \\
dest13    & suppl18    & plant9                &  3.22581 \\
\hline
dest14         &  suppl19  & plant5  &   1.35355 \\
dest14         &  suppl20  & plant5 &   3.22581 \\
\hline
dest15 & suppl7 & plant1 & 0.00433257 \\
\end{tabular}
}
\end{table}


\begin{table}[ht!]
\centering
{\tiny
\caption{Optimal solution of  RO-ell  with $\Omega=2.75$ and  all continuous variables. The table shows the optimal number of booked vehicles $x_{ijk}$ from plant $i\in\mathscr{O}_k$ of supplier $k$ to destination $j\in\mathscr{D}$.}
 \label{tabROell}
\begin{tabular}[l]{@{}llll}
\textbf{Destination} $j\in\mathscr{D}$ & \textbf{Supplier} $k\in\mathscr{K}$ & \textbf{Plant} $i\in\mathscr{O}_k$ & $x_{ijk}$\\
\hline
dest1 & suppl1               &  plant3   &                 9.27 \\
\hline
dest2    & suppl1               & plant3      &              38.05 \\
dest2    & suppl4        & plant4 & 2.13  \\
dest2    & suppl6 & plant10 & 10.52\\
dest2    & suppl12                & plant14      &              3.23 \\
dest2    & suppl16       & plant16 &         3.23    \\
\hline
dest3         & suppl3              & plant8          &   5.88      \\
dest3         & suppl6               & plant10          &   1.64 \\
dest3         & suppl11               & plant13              & 2.16\\
dest3         & suppl15              & plant8          &   2.72 \\
dest3         & suppl17            & plant17           &   3.23 \\
dest3         & suppl19            & plant5         &   0.84 \\
dest3         & suppl20               & plant5         &   0.84 \\
\hline
dest4 & suppl13 & plant15  &  3.23  \\
\hline
dest5 & suppl3 & plant8 &  0.57 \\
dest5 & suppl15 & plant8 &  0.5\\
\hline
dest6         & suppl10               & plant12          &   1.83 \\
dest6         & suppl14             & plant12          &   1.97 \\
dest6         & suppl24            & plant12          &   1.96 \\
\hline
dest7 & suppl19   &  plant5 &    2.27 \\
dest7 & suppl20   & plant5&    2.27 \\
\hline
dest8   & suppl5   & plant9    &              1.96 \\
dest8   & suppl18               & plant9   &              1.26  \\
dest8   & suppl9               &  plant11    &              4.84 \\
dest8   & suppl22                & plant19    &              3.23 \\
\hline
dest9 & suppl2 & plant7 &    2.89   \\
dest9 & suppl10  & plant12 &   1.21\\
dest9 & suppl14    & plant12 &   1.26 \\
dest9 & suppl21  & plant18 &   3.23\\
dest9 & suppl23  & plant7 &   0.86\\
dest9 & suppl24   & plant12 &   1.26 \\
\hline
dest10      & suppl1               & plant3                  &   1.06 \\
\hline
dest11           & suppl11 & plant13   &  1.68\\
\hline
dest13    & suppl2    & plant7                &  3.56 \\
dest13    & suppl5    & plant9                &  4.49 \\
dest13    & suppl18    & plant9                &  1.97 \\
dest13    & suppl23    & plant7                &  0.82 \\
\hline
dest14         &  suppl19  & plant5 &   0.11 \\
dest14         &  suppl20  & plant5  &   0.11 \\
\hline
dest15       & suppl7  & plant1 &  3.23 \\
dest15       & suppl8   & plant2 &    3.21 \\
dest15       & suppl11  & plant13 &    5.19 \\
\hline
\end{tabular}
}
\end{table}

\begin{figure}[ht!]
\centering
\includegraphics[width=\textwidth]{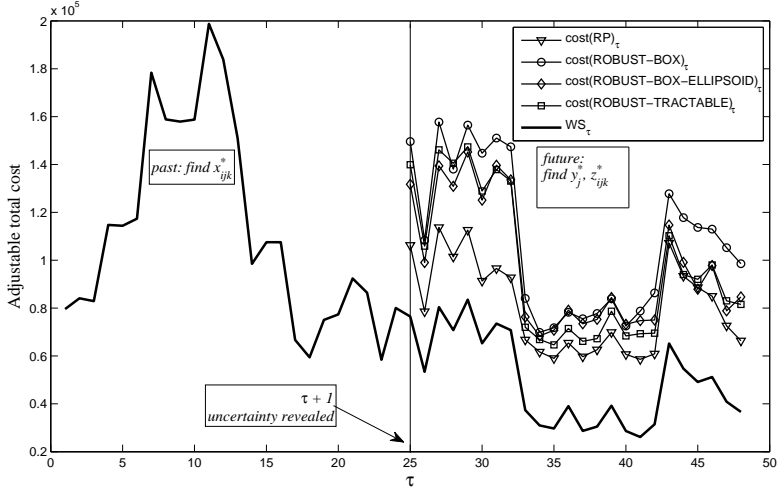}
\vskip3truemm\noindent
\caption{Adjustable total costs $\mathbf{{cost}}_{\mathbf{1},\tau},\dots,\mathbf{{cost}}_{\mathbf{4},\tau}$ obtained by solving model (\ref{o1})-(\ref{i7}) with all continuous variables using the nonadjustable decision variables $x_{ijk}^{\ast}$  given by methods $\mathbf{M_1},\dots,\mathbf{M_4}$, for an increasing value of $\tau=24,\dots,47$.
Results are compared with the cost of the  problem where a full information on the realization of vectors $b$ and $d$ is available ($\mathbf{{cost}}_{\tau}=WS_{\tau}$). }
\label{resultsadjustablefig}
\end{figure}

\begin{table}[ht!]
\centering
{\tiny
\caption{Optimal value of the objective function (\ref{o1})  with the adjustable strategy where  method $\mathbf{M_1},\dots,\mathbf{M_5}$ is used to determine the nonadjustable variables $x_{ijk}^{\ast}$.
Last column ($\mathbf{{cost}}_{\tau}$) refers to the cost of problem  where a full information on the realization of vectors $b$ and $d$ is available. The aggregated cost $\mathbf{{cost}_{m}}$, $\mathbf{m}=\mathbf{1},\dots,\mathbf{5}$,  is reported on the last line. All variables are supposed continuous.}
 \label{tabadjcosts}
\begin{tabular}[l]{@{}lrrrrrr}
$\tau$ & $\mathbf{{cost}}_{\mathbf{1},\tau}$ & $\mathbf{{cost}}_{\mathbf{2},\tau}$ & $\mathbf{{cost}}_{\mathbf{3},\tau}$ & $\mathbf{{cost}}_{\mathbf{4},\tau}$ & $\mathbf{{cost}}_{\mathbf{5},\tau}$ & $\mathbf{{cost}}_{\tau}$ \\
\hline
24 	& 106268 &	149596 &	131733 &	139920 &	$\infty$ &	76464\\
25 	& 78537 &	108328 &	99046 &	105893 &	$\infty$	& 53386\\
26 	& 113678	& 157730 &	139411 &	146097	& $\infty$	& 80399\\
27 & 101405	& 138061	& 130867 & 	140536	& $\infty$ &	70792\\
28 	& 112587& 156452	& 145235 &	147368	& $\infty$	& 83548\\
29 	  &  91334 &	144716 &	125087 &	128931 &	$\infty$	& 65283\\
30 	& 96586	& 151069	& 139699	 & 137865	& $\infty$	& 73541\\
31 	& 92735	& 147384	& 133558	& 133085	& $\infty$	& 70808\\
32 	& 66738	& 84055	& 76295	& 72009	& $\infty$	& 37361\\
33 & 61685	& 69863	& 68364	& 66863 &	$\infty$	& 30987\\
34 & 	58994 & 	71826	& 70761 & 	64654	& $\infty$	& 29690\\
35 	& 65425	& 78257	& 79132	& 71442	& $\infty$	 & 39029\\
36 	& 59721& 75545	& 73450	& 66142	& $\infty$	& 28702\\
37 	& 62597	& 77775	& 75307 & 67223	& $\infty$	& 30513\\
38  & 	69873	&  83871 & 	84397 & 	78665	& $\infty$ & 39240\\
39 	& 60691	&  72534	& 73205	& 68384	& $\infty$ & 28646\\
40 	& 58604	&  78762 & 	74779	& 69311	& $\infty$ & 26157\\
41 	& 60879	& 86342	& 75022	& 69542	& $\infty$ & 31475\\
42 	& 107178	& 127786 & 	114733	& 110211	& $\infty$ & 65190\\
43  & 	93369	& 117811	& 99104	& 94021	& $\infty$ & 54748\\
44  & 	88607	& 113692	& 87892	& 92033	& $\infty$ & 49120\\
45 	 & 84985	& 113013	& 97798	& 98013	& $\infty$ & 51223\\
46  & 	72613	& 105263	& 78852	& 83042	& $\infty$ & 40894\\
47  & 	66353	& 98546	& 84675	& 81543	& $\infty$ & 36623\\
\hline
 & $\mathbf{cost_{1}}$ & $\mathbf{cost_{2}}$ & $\mathbf{cost_{3}}$ & $\mathbf{cost_{4}}$ & $\mathbf{cost_{5}}$ & $\mathbf{cost}$\\
\hline
 & 1931443 &	2608276 &	2358398 &	2332794 & $\infty$ & 1193820
 \end{tabular}
}
\end{table}


In order to make a fair comparison with the SP methodology,  we compute total costs via the scenario based  framework proposed in Section \ref{sec:adjust}.
Our methodology allows us to understand  the cost saving of the SP approach when compared to  RO, and to quantify the value of a more conservative strategy which does not require a negotiation  with third-party providers  every week, but  keeps the same solution strategy for longer periods.
For this purpose, we consider a subset $\bar{\mathscr{S}}=\left\{1,\dots,24\right\}\subset \mathscr{S}$ of scenarios from the 48 historical data.
For each $\tau=24,\dots,47$ and for each method $\mathbf{M_1},\dots,\mathbf{M_5}$, we compute nonadjustable decision variables $x_{ijk}^{\ast}$ using only the information contained in $\hat{d}^1,\hat{d}^2,\ldots,\hat{d}^{\tau},\quad \hat{b}^1,\hat{b}^2,\ldots,\hat{b}^{\tau}$.
When vectors $\hat{d}^{\tau + 1},\hat{b}^{\tau+1}$ become available, we solve model (\ref{o1})-(\ref{i7}) fixing the nonadjustable decision variables $x_{ijk}^{\ast}$  just obtained, allowing to determine the adjustable decision variables $y^{\ast}_j,z_{ijk}^{\ast}$. The adjustable total costs $\mathbf{{cost}}_{\mathbf{1},\tau},\dots,\mathbf{{cost}}_{\mathbf{5,\tau}}$ and $\mathbf{{cost}}_{\tau}$  are shown in Figure \ref{resultsadjustablefig} and results are listed in Table \ref{tabadjcosts}.
Last line reports the aggregated cost $\mathbf{{cost}_{m}}$, $\mathbf{m}=\mathbf{1},\dots,\mathbf{5}$ over 24 weeks.
The results provide an important information to the firm about the cost saving in case of SP or RO solution procedures are implemented in practice over 6 months:
SP approach allows a saving of $35.04\%$	compared to RO-box, $22.10\%$	compared to RO-ell ($\Omega=2.75$) and $20.77\%$ compared to the computationally tractable problem trSOCP.
Nevertheless, RO solutions are immunized over all realization of uncertain parameters allowing to the company to keep the same contract for longer periods without the necessity of a weekly negotiation and an adjustment of the plan when the booked vehicles are not sufficient. Only in case the observations of demands and costs in $\tau + 1$ are worse than their history up to $\tau$, should the RO solution strategy also be renegotiated.
However this would be not the case when $(\hat{d}^{\tau + 1},\hat{b}^{\tau+1})$ corresponds to an extremely bad scenario:
the results are shown in Table \ref{tabextreme} where we can see the better performance of RO approaches with respect to SP allowing a saving of $3.32\%$.

On the other hand, the  adjustable RO approach $\mathbf{{M}_{5}}$ is no longer implementable since $\varphi^\tau\left(\lambda\left(\bar{d}^{\tau+1}\right)\right)>0$ and consequently $\mathbf{{cost}}_{\mathbf{5},\tau}=\infty$ for $\tau=24,\dots,47$.

As expected the lowest cost  is given by the  wait-and see WS problem allowing a saving of $38.19\%$ compared to SP, since a full information on the realization of vectors $b$ and $d$ is available.

\begin{table}[ht!]
\begin{center}
\caption{Optimal value of the objective function (\ref{o1})  with the adjustable strategy where  method $\mathbf{M_1},\dots,\mathbf{M_5}$ is used to determine the nonadjustable variables $x_{ijk}^{\ast}$  in case of worst scenario $(\hat{d}^{\tau + 1},\hat{b}^{\tau+1})=(\bar{d}_j + \gamma_j \cdot \bar{d}_j,\bar{b}_j + \sigma \cdot \bar{d}_j)$, $j\in\mathscr{D}$.
Last column ($\mathbf{{cost}}_{\tau}$) refers to the cost of the problem  where a full information on the realization of vectors $b$ and $d$ is available.}
\label{tabextreme}
\begin{tabular}[l]{@{}l|llll|l}
SP & & RO & & & WS \\
 $\mathbf{cost_{1}}$ & $\mathbf{cost_{2}}$ & $\mathbf{cost_{3}}$ & $\mathbf{cost_{4}}$ & $\mathbf{cost_{5}}$ & $\mathbf{cost}$\\
\hline
350430.70 &	339150.92 &	346527.34 &	352759.79 & $\infty$ & 198756.53
 \end{tabular}
\end{center}
\end{table}

Total CPU time, in seconds, spent in solving the optimization problems $\mathbf{M_m}$,  $\mathbf{m}=\mathbf{1},\dots,\mathbf{5}$ and (\ref{o1})-(\ref{i7}) including the possible infeasibility detection of the latter, are reported  in Table  \ref{tabadjCPU}. Results  show  the higher computational complexity of the SP approach $\mathbf{M_1}$ and the adjustable approcahes  $\mathbf{M_4}$ and $\mathbf{M_5}$ as compared  to the robust box constrained problem $\mathbf{M_2}$ or the robust box-ellipsoidal constrained problem $\mathbf{M_3}$.

\begin{table}[ht!]
\centering
{\tiny
\caption{Total CPU time, in seconds, spent in solving the optimization problems $\mathbf{M_m}$,  $\mathbf{m}=\mathbf{1},\dots,\mathbf{5}$. Last column ($\mathbf{{cost}}_{\tau}$) refers to the CPU time of problem where a full information on the realization of vectors $b$ and $d$ is available.}
 \label{tabadjCPU}
\begin{tabular}[l]{@{}lllllll}
$\tau$ & $\mathbf{{CPU}}_{\mathbf{1},\tau}$ & $\mathbf{{CPU}}_{\mathbf{2},\tau}$ & $\mathbf{{CPU}}_{\mathbf{3},\tau}$ & $\mathbf{{CPU}}_{\mathbf{4},\tau}$ & $\mathbf{{CPU}}_{\mathbf{5},\tau}$  & $\mathbf{{CPU}}_{\tau}$ \\
\hline
24 	& 8.8356	  & 0.0312	& 0.2156	& 32.0156 	& 32.0156 &	0.0156 \\
25 	& 9.1856	  & 0.0312	& 0.1856	& 29.8856	& 29.8856&	0.0156 \\
26 	& 8.1656	  & 0.0312	& 0.2156	& 30.0156	& 30.0156 &	0.0156 \\
27 	& 9.8856	  & 0.0312	& 0.2256	& 32.9856	& 32.9856 &	0.0156 \\
28 	& 10.0156	& 0.0312	& 0.2256	& 77.2756	&  77.2756 &	0.0156 \\
29	& 10.9756	& 0.0312	& 0.1956	& 83.5056	& 83.5056 &	0.0156 \\
30	& 11.8956	& 0.0312	& 0.2156	& 31.1456	& 31.1456 &	0.0156 \\
31 	& 11.1456	& 0.0312	& 0.1956	& 59.3956	& 59.3956 &	0.0156 \\
32 	& 10.7456	& 0.0312	& 0.2156	& 73.3956	& 73.3956 &	0.0156 \\
33 	& 12.1756	& 0.0312	& 0.1956	& 71.2056	& 71.2056 &	0.0156 \\
34 	& 11.4256	& 0.0312	& 0.2156	& 76.5256	& 76.5256 &	0.0156 \\
35 	& 10.8256	& 0.0312	& 0.2156	& 65.4056	& 65.4056 &	0.0156 \\
36 	& 11.5856	& 0.0312	& 0.1856	& 58.2256	&  58.2256 &	0.0156 \\
37 	& 12.3156	& 0.0312	& 0.1956	& 91.2756	& 91.2756 &	0.0156 \\
38 	& 11.3956	& 0.0312	& 0.2256	& 53.5856	& 53.5856 &	0.0156 \\
39 	& 10.7256	& 0.0312	& 0.2156	& 83.5456	& 83.5456 &	0.0156 \\
40 	& 11.1056	& 0.0312	& 0.2156	& 39.6156	& 39.6156 &	0.0156 \\
41 	& 10.9856	& 0.0312	& 0.2026	& 65.8256	& 65.8256 &	0.0156 \\
42 	& 10.9656	& 0.0312	& 0.2156	& 62.3156	& 62.3156 &	0.0156 \\
43	& 11.2756	& 0.0312	& 0.2156	& 89.1056	& 89.1056 &	0.0156 \\
44 	& 11.2856	& 0.0312	& 0.2026	& 75.4256	& 75.4256 &	0.0156 \\
45 	& 11.6056	& 0.0312	& 0.2156	& 89.4656	& 89.4656 &	0.0156 \\
46 	& 11.3356	& 0.0312	& 0.2026	& 73.7656	& 73.765 &	0.0156 \\
47 & 14.0056	& 0.0312	& 0.2156	& 104.4256	& 104.4256 &	0.0156 \\
 \end{tabular}
}
\end{table}



\begin{figure}[ht!]
\centering
\includegraphics[width=\textwidth]{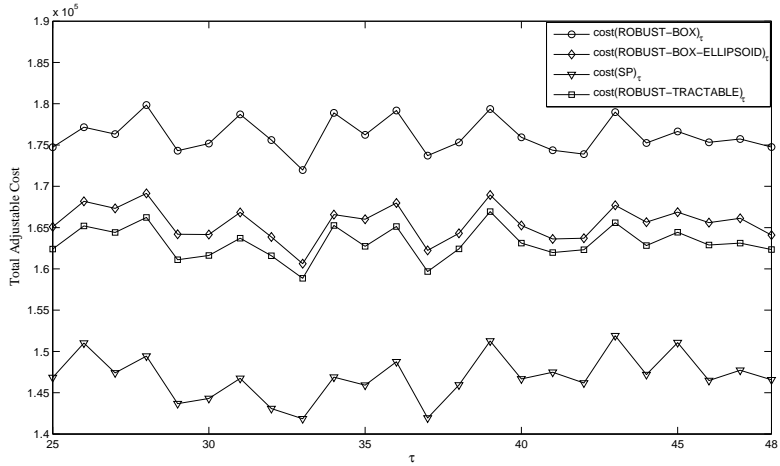}
\vskip3truemm\noindent
\caption{Adjustable total costs $\mathbf{{cost}}_{\mathbf{1},\tau},\dots,\mathbf{{cost}}_{\mathbf{4},\tau}$ obtained by solving model (\ref{o1})-(\ref{i7}) using the nonadjustable decision variables $x_{ijk}^{\ast}$  respectively given by methods $\mathbf{M_1},\dots,\mathbf{M_4}$, for an increasing value of $\tau=24,\dots,47$. Results are averaged over 2000 simulations randomly generated by  a Monte Carlo procedure.}
\label{resultsadjustableA}
\end{figure}

We finally validate the performance of the approaches $\mathbf{M_1},\dots,\mathbf{M_5}$ over 2000 simulations randomly generated by  a Monte Carlo procedure. Uncertain buying cost values are obtained by sampling from a uniform distribution  in the interval $\left[\bar{b}_j - \sigma \cdot \bar{b}_j ,\bar{b}_j + \sigma \cdot \bar{b}_j \right],\ j\in\mathscr{D}$, with a  given deviation level of $\sigma=20\%$. Uncertain demand values are obtained by sampling from a uniform distribution  in the interval $\left[\bar{d}_j - \gamma_j \cdot \bar{d}_j ,\bar{d}_j + \gamma_j \cdot \bar{d}_j \right],\ j\in\mathscr{D}$
with  $\gamma_j=\max_{s=1,\dots,48} \hat{d}_j^s -\bar{d}_j, \; j\in\mathscr{D}$ (see Table \ref{tab2}).
Total costs over 24 weeks over 2000 simulations are reported in Table \ref{tabadjSIMULATION} and shown in Figure \ref{resultsadjustableA}. The numerical results show that the SP approach allows a saving of $19.88\%$	compared to RO-box, $12.69\%$	compared to RO-ell and $11.07\%$ compared to computationally tractable robust formulation trSOCP. Again, the adjustable method $\mathbf{M_5}$ is not implementable.

\begin{table}[ht!]
\begin{center}
\caption{Total cost of optimization problems $\mathbf{M_m}$,  $\mathbf{m}=\mathbf{1},\dots,\mathbf{5}$ and (\ref{o1})-(\ref{i7}) over 24 weeks over 2000 simulations randomly generated by  a Monte Carlo procedure. }
 \label{tabadjSIMULATION}
\begin{tabular}[l]{@{}l|llll}
SP & & RO & &  \\
 $\mathbf{cost_{1}}$ & $\mathbf{cost_{2}}$ & $\mathbf{cost_{3}}$ & $\mathbf{cost_{4}}$ & $\mathbf{cost_{5}}$ \\
\hline
	3526310.409	& 4227395.303	 & 3974122.436	& 3915887.398 & $\infty$
 \end{tabular}
\end{center}
\end{table}


\section{Conclusions}
\label{sec:Conclusions}
In this paper we have analyzed the effect of two modelling approaches, stochastic programming (SP) and robust optimization (RO) for a  supply planning problem under uncertainty.
The problem has been formulated and  solved both via a two-stage stochastic programming and robust optimization  models with different uncertainty sets. 

The goal of SP is to compute the minimum expected cost based on the specific probability distribution of the uncertain parameters based on a  set of  scenarios.

For RO we have firstly considered  static approaches with random parameters belonging to box or ellipsoidal uncertainty sets in compliance with the data used to generate scenarios for SP, and secondly dynamic approaches, via the concept of adjustable robust counterpart ARC.

The choice of the box uncertainty set is preferable only if the feasibility of all the constraints is highly required, but this certainty of constraint satisfaction results in a higher cost. A  less conservative outcome has been obtained with  a box-ellipsoidal uncertainty set that  requires the solution of a second-order cone program SOCP.
The main advantage of the considered RO formulations, is that they can be solved in polynomial time and theoretical guarantees for the quality of the solution are provided, which is not the case with the aforementioned SP formulation.


A  scenario based framework methodology for a fair comparison between SP and RO has been proposed, which can be applied to any optimization problem affected by uncertainty.

The efficiency of the methodology has been illustrated for a supply planning problem to optimize
vehicle-renting and transportation activities involving uncertainty on demands and buying costs for
extra-vehicles.
The methodology allows to understand what is the cost saving of the SP approach when compared to the RO approach and to quantify the value of a more conservative strategy which does not require a  negotiation  with suppliers or third-party providers  every week.

\section{Annex}
\label{sec:Annex}
\begin{table}[ht!]
\centering
{\tiny
\caption{Set of suppliers $\mathscr{K}$ and set of their plants ${O}_k, \ k \in\mathscr{K}$.}
 \label{tab1}
\begin{tabular}[l]{@{}lll}
& \textbf{Supplier} $k\in\mathscr{K}$ & \textbf{Plant} $i\in\mathscr{O}_k$\\
\hline
 & suppl1 &	plant1 \\
& 	& plant2\\
&	& plant3 \\
&	& plant4\\	
& & plant5\\
&	& plant6\\
	\hline
 & suppl2	& plant7 \\
\hline
 & suppl3	& plant8\\
\hline
 & suppl4	& plant4\\
\hline
 & suppl5 &	plant9 \\
\hline
 & suppl6	& plant10 \\
	&        &   plant6 \\
\hline
 & suppl7	& plant1 \\
\hline
 & suppl8	& plant1 \\
	 &                             & plant2 \\
\hline
 & suppl9 &	plant11 \\
\hline
 & suppl10	& plant12  \\
\hline
 & suppl11	& plant13 \\
\hline
 & suppl12	& plant14 \\
\hline
 & suppl13	& plant15 \\
\hline
 & suppl14	& plant12 \\
\hline
 & suppl15	 & plant8 \\
\hline
 & suppl16 &	plant16 \\
\hline
 & suppl17 &	plant17 \\
\hline
 & suppl18	& plant9 \\
\hline
 & suppl19	& plant5 \\
&	& plant15\\
\hline
 & suppl20	& plant5 \\
\hline
 & suppl21 &	plant18 \\
\hline
 & suppl22	& plant19 \\
\hline
 & suppl23	& plant7 \\
\hline
 & suppl24	& plant12 \\
 \end{tabular}
}
\end{table}

\begin{table}[ht!]
\centering
{\tiny
\caption{List of destinations (cement factories) with relative expected emergency costs $\bar{b}_j$, maximum capacities which can be booked $g_j$, $j \in\mathscr{D}$ and values of deviation for defining  the box constraints for demand (\ref{box1})  and buying cost (\ref{box2}).}
 \label{tab2}
\begin{tabular}[l]{@{}llllll}
& Destination $j\in\mathscr{D}$ & Expected emergency cost $\bar{b}_j$ & maximum booking capacity $g_j$ & $\rho_1 G_j$ & $\rho_2 F_j$\\
 \hline
 & dest1 & 72.61 & 422.95 & 226.28 & 9.70\\
 & dest2 &  70.58 & 2054.55 & 1237.99 & 11.17\\
 & dest3 & 68.01 & 1330.67 & 428.53 & 10.75\\
 & dest4 & 64.94 & 453.64 & 241.69 & 10.19\\
 & dest5 & 73.52  & 613.41 & 92.40 & 11.28\\
 & dest6 & 58.57 & 695.24 & 264.61 & 10.96 \\
 & dest7 & 69.83  & 443.14 & 174.08 & 11.65\\
 & dest8& 66.32 & 815.36 & 312.29 &  10.32 \\
 & dest9 & 62.63 & 933.33 & 397.64 & 10.65 \\
 & dest10 & 68.22 & 319.79 & 175.75 & 9.94 \\
 & dest11 & 48.92 & 443.11 & 181.37 & 12.39\\
 & dest12 & 50.04  & 760.11 & 326.82 & 10.03 \\
 & dest13 & 73.07 & 381.20 & 228.38 & 9.79\\
 & dest14 & 59.93 & 498.33 & 191.78 & 10.97 \\
 & dest15 & 55.63  & 232411.75 & 564.67 & 10.94\\
 \end{tabular}
}
 \end{table}
\begin{table}[ht!]
\centering
{\tiny
\caption{Minimum  $r_k$ and average maximum $v_k$ requirement capacity of supplier $k \in\mathscr{K}$.}
 \label{tabupper}
\begin{tabular}[l]{@{}llll}
& Supplier $k\in\mathscr{K}$ & $r_k$ & $v_k$  \\
\hline
 & suppl1  & 1057.69 & 1500  	 \\
 & suppl2 & 0	& 200  \\
 & suppl3 & 0	& 200 \\
 & suppl4 &	0 & 66 \\
 & suppl5 & 0 & 200 	 \\
 & suppl6	& 0 & 376.92  \\
 & suppl7 &	0 & 100  \\
 & suppl8	&  0 & 100 \\
 & suppl9 & 0 & 150	 \\
 & suppl10 &	0 & 94   \\
 & suppl11 & 0	& 280  \\
 & suppl12 & 0	& 100  \\
 & suppl13 & 0	&  100  \\
 & suppl14 & 0	& 100  \\
 & suppl15	& 0 & 100  \\
 & suppl16  & 0 & 100 	 \\
 & suppl17 & 0 &	100   \\
 & suppl18	& 0 & 100   \\
 & suppl19	& 0 & 100  \\
 & suppl20	 & 0 &  100  \\
 & suppl21 & 0 & 100 	\\
 & suppl22	& 0 & 100   \\
 & suppl23	& 0 &   50 \\
 & suppl24	&  0 &  100  \\
 \end{tabular}
}
 \end{table}

 \section*{Acknowledgements}
The authors would like to thank the \emph{Italcementi} Logistics Group, in particular Dott.\ Luca Basaglia, Francesco Bertani, and Flavio Gervasoni  for the description of the problem and the historical data provided.
F. Maggioni and M. Bertocchi acknowledge the 2014-2015 University of Bergamo grants.
The work of Florian Potra was supported in part by the National Science Foundation under Grant No. DMS-1311923.



\begin{thebibliography}{100}
\bibliographystyle{ormsv080}

\bibitem{AK05}  Alidaee, B., \& Kochenberger, G.A. (2005)  A note on a simple dynamic programming approach to the single-sink, fixed-charge transportation problem. {\em Transportation Science}, 39(1), 140-143.
\bibitem{BTN99}  Ben-tal, A., \& Nemirovski, A. (1999) Robust solutions of uncertain linear programs. \emph{	Operations Research Letters}, 25, 1-13.
\bibitem{BGGN04}  Ben-Tal, A.,  Goryashko, A.,  Guslitzer, E., \&    Nemirovski,  A. (2004) Adjustable robust solutions of uncertain linear programs. \emph{Mathematical Programming}, 99(2), 351-376.
\bibitem{BTN00} Ben-Tal, A., \&  Nemirovski, A. (2000) Robust solutions of linear programming problems contaminated with uncertain data. \emph{Mathematical  Programming (Series B)} 88, 411-424.
\bibitem{BTEGN09} Ben-Tal, A., El-Ghaoui, L., \&  Nemirovski, A. (2009) \emph{Robust optimization}. Princeton University Press ISBN  978-0-691-14368-2.
\bibitem{BS04} Bertsimas, D., \&   Sim, M. (2004)   The price of robustness. \textit{Operations Research}, 52(1), 35-53.
\bibitem{BG12} Bertsimas, D., \&  Goyal, V. (2012)  On the power and limitations of affine policies in two-stage adaptive optimization. \textit{Mathematical Programming (Series A)},  134, 491-531.
\bibitem{BBC2011} Bertsimas, D., Brown, D.B., \& Caramanins, C. (2011) Theory and applications of robust optimization. \textit{Siam Review}, 53(3), 464-501.
\bibitem{BB2010} Bertsimas, D., \& Caramanis,  C. (2010) Finite adaptibility in multistage linear optimization. \emph{IEEE Transactions
on Automatic Control}, 55(12), 2751–2766.
\bibitem{CSS2014} Chen, X.,  Shum, S., \&  Simchi-Levi D. (2014), Stable and Coordinating Contracts for a Supply Chain with Multiple Risk-Averse Suppliers. \textit{Production and Operations Management}.  23(3),  379-392.

\bibitem{BG2013} Bertsimas, D., \& Georghiou, A. (2015) Design of near optimal decision rules in multistage adaptive mixed-integer
optimization.  \textit{Operations Research}, 63(3), 610-627.
\bibitem{BG2014} Bertsimas, D.,  \& Georghiou, A. (2014) Binary decision rules for multistage adaptive mixed-integer optimization.
 \emph{Optimization Online}.
\bibitem{BG2010} Bertsimas, D.,  \& Goyal, V. (2010) On the power of robust solutions in two-stage stochastic and adaptive
optimization problems. \emph{Mathematics of Operations Research}, 35(2), 284–305
\bibitem{BL11} Birge, J.R., \&  Louveaux, F. (2011)  \emph{Introduction to stochastic programming}, Springer-Verlag, New York.
\bibitem{CP96} Cheung, R.K., \& Powell, W.B. (1996) Models and Algorithms for Distribution Problems with Uncertain Demands. {\em Transportation Science}, 30, 43-59.
\bibitem{CL77} Cooper, L.,  \&   LeBlanc, L.J. (1997)  Stochastic transportation problems and other network related convex problems. {\em Naval Research Logistics Quarterly}, 24, 327-336.
\bibitem{1} Coyle, J.J., Bardi, E.J., \&  Langley, C.J. (2003) \emph{The Management of Business Logistics—A Supply Chain Perspective}. South- Western Publishing, Mason.
\bibitem{2} Crainic, T.G., \&   Laporte, G., 1997, Planning models for freight transoportation. \emph{European Journal of Operational Research}, 97, 409-438.

\bibitem{2a} Crainic, T.G., Gobbato, L., Perboli, G., \&  Rei, W. (2014) Logistics capacity planning: a stochastic bin packing formulation and a progressive hedging heta-heuristic, \emph{CIRRELT} 2014-66.

\bibitem{CL97}  Crainic, T.G., \&     Laporte, G. (1997) Planning models for freight transportation. {\em European Journal of Operational Research},  97(3), 409-438.
\bibitem{D98}   Dupacova, J.  (1998) Reflections on robust optimization.  Marti K,  Kall P, eds. {\em Stochastic Programming Methods and Technical Applications}, (Springer Verlag, Berlin Heidelberg),  111-127.
\bibitem{EL97} 	 El Ghaoui L.,	\& Lebret, H. (1997) Robust Solutions to Least-Squares Problems with Uncertain Data. \emph{SIAM Journal on Matrix Analysis and Applications archive}, 18(4),  1035-1064
\bibitem{EOL98} 	 El Ghaoui, L., Oustry, F., \&	Lebret, H. (1998) Robust Solutions to uncertain semidefinite programs. \emph{SIAM Journal of Optimization}, 9, 33-52.
\bibitem{G93}  Graves, S.C.,   Rinnooy Kan, A.H.G., \&  Zipkin, P.H. (1993) {\em Logistics of Production and Inventory},  Handbooks in Operations Research and Management Science, 4,    ISBN: 978-0-444-87472-6
\bibitem{KW} Kaut, M., \& Wallace, S.W. (2007) Evaluation of scenario generation methods for stochastic programming. \emph{Pacific Journal of Optimization}, 3(2), 257–271.
\bibitem{3}Knemeyer, A.M., Corsi, T.M., \& Murphy, P.R. (2003) Logistics outsourcing relationships: Customer perspectives. \emph{Journal of Business Logistics}, 24(1), 77–109.
\bibitem{KWG11} Kuhn, D.,  Wiesemann, W.,  \& Georghiou, A. (2011) Primal and dual linear decision rules in stochastic and robust optimization. {\em Mathematical Programming},  130(1), 177-209.
\bibitem{HKW2015}  Hanasusanto, G.A., Kuhn, D., \& Wiesemann, W. (2015) K-Adaptability in Two-Stage Robust Binary Programming. \textit{Operations Research}, 63(4), 877 - 891{\em Optimization Online}.
\bibitem{LW97} Lamar, B.W., \&  Wallace, C.A.  (1997) Revised-modified penalties for fixed charge transportation problems. {\em Management Science}, 43(10), 1431-1436.
\bibitem{LSP90} Lamar, B.W.,  Sheffi, Y., \&  Powell, W.B. (1990)  A capacity improvement lower bound for fixed charge network design problems. {\em Operations Research}, 38(4), 704-710.
\bibitem{4}Lieb, R.C., \& Bentz, B.A. (2004) The use of third-party logistics services by large American  manufacturers: The 2003 survey. \emph{Transportation Journal}, 43(3), 24–33.
\bibitem{5}Lieb, R.C., \& Bentz, B.A. (2005) The use of third-party logistics services by large American manufacturers: The 2004 survey. \emph{Transportation Journal}, 2, 5–15.
\bibitem{6}Lieb, R.C., \& Miller, J. (2002) The use of third-party logistics services by large US manufacturers, the 2000 survey. \emph{International Journal of Logistics: Research and Applications}, 5(1), 1–12.
\bibitem{7}Lieb, R.C. (1992) The use of third-party logistics services by large American manufacturers. \emph{Journal of Business Logistics}, 13(2), 29–42.
\bibitem{8} Lieb, R.C., \& Randall, H.L. (1999) Use of third-party logistics services by large US manufacturers in 1997 and comparisons with previous years. \emph{Transport Reviews}, 19(2), 103–115.


\bibitem{MKB} Maggioni, F.,  Kaut, M., \&  Bertazzi, L. (2009) Stochastic Optimization models for a single-sink transportation problem. {\em Computational Management Science}, 6, 251-267.
\bibitem{MAB2014} Maggioni, F.,  Allevi, E., \&  Bertocchi, M. (2014)  Bounds in Multistage Linear stochastic Programming. {\em Journal of Optimization, Theory and Applications}, 163(1), 200-229.
\bibitem{9}Maltz, A., \& Ellram, L.M. (1997) Total cost of relationship: An analytical framework for the logistics outsourcing decision. \emph{Journal of Business Logistics}, 18(1), 45–66.
\bibitem{10}Marasco, A. (2008) Third-party logistics: A literature review. \emph{International Journal of Production Economics}, 113, 127-147.
\bibitem{PoTo2003}  Powell, W.B., \&    Topaloglu, H.  (2003) Stochastic Programming in Transportation and Logistics, in \textit{Handbooks in Operations Research and Management Science}, 10, 555-635.
\bibitem{RS03} Ruszczinski,  A., \&  Shapiro, A. {\em Stochastic Programming}, Elsevier, Amsterdam (2003).
\bibitem{12} Razzaque, M.A., \& Sheng, C.C.  (1998) Outsourcing of logistics functions: A literature survey. \emph{International Journal of Physical Distribution \& Logistics Management}, 28(2), 89–107.
\bibitem{13}Simchi-Levi, D., \&  Simchi-Levi, E. (2003) Inventory Positioning: Exploring Push and Pull Supply Chains.  \emph{Parcel Shipping \& Distribution}.
\bibitem{14} Simchi-Levi, D., \&  Simchi-Levi, E. (2004) Inventory Optimization: The Last Frontier. \emph{Inbound Logistics}.
\bibitem{15} Simchi-Levi, D. (2014), OM Research: From Problem Driven to Data Driven Research. \emph{Manufacturing and Service Operations Management}, 16(1),  2-10.
\bibitem{SCL2011} Solyali, O.,   Cordeau, J.F., \&  Laporte, G. (2001)  Robust Inventory Routing Under Demand Uncertainty. \emph{Transportation Science}, 46(3), 327-340.
\bibitem{So73}  Soyster, A.L. (1973)  Convex programming with set-inclusive constraints and applications to inexact linear programming. {\em Operations Research}, 43(2), 264-281.
\bibitem{LV02}  Van Landeghem,  H., \&  Vanmaele, H. (2002) Robust planning: a new paradigm for demand chain planning. {\em Journal of Operations Management}, 20(6), 769-783.
\bibitem{YL00} Yu, C., \&  Li, H. (2000) A robust optimization model for stochastic logistic problems. {\em International Journal of Production Economics}, 64(1-3), 385-397.





\end{thebibliography}
\end{document}